\documentclass[12pt]{article}
\usepackage{amssymb,latexsym,amsmath}
\def\C{\mathbb C}
\def\R{\mathbb R}
\def\N{\mathbb N}
\def\Z{\mathbb Z}

\numberwithin{equation}{section}

\newtheorem{thm}{Theorem}[section]
\newtheorem{lem}{Lemma}[section]

\newtheorem{definition}{Definition}[section]
\newtheorem{prop}{Proposition}[section]

\begin{document}
\title{Non-real zeros of linear differential polynomials}
\author{J K Langley}
\maketitle
\sffamily
\begin{abstract}
Let $f$ be a real entire function with finitely many non-real zeros,
not of the form $f = Ph$ with $P$ a polynomial and $h$ in the
Laguerre-P\'olya class.
Lower bounds are given for the number of non-real zeros of
$f'' + \omega f$, where $\omega$ is a positive real constant.
\end{abstract}

\section{Introduction}

This paper concerns non-real zeros of linear differential polynomials
in real entire functions with real zeros. 
For each
non-negative integer $p$ the class $V_{2p}$ 
\cite{HelW1,HelW2,SS} consists of all entire functions
$$
f(z) = g(z) \exp( -a z^{2p+2} ) ,
$$
where $a \geq 0$ is real and $g$ is a real entire function 
with real zeros of 
genus at most $2p+1$ \cite[p.29]{Hay2}. The classes $U_{2p}, p \geq 0$, are
then given by $U_0 = V_0$ and
$U_{2p} = V_{2p} \setminus V_{2p-2}$ for $ p \geq 1$.
Moreover, $U_{2p}^*$ is the  class of entire functions $f= Ph$, 
where $h \in U_{2p} $ and $P$ is a real polynomial without
real zeros \cite{EdwH}, so that every real
entire function of finite order with finitely many non-real zeros
belongs to $U_{2p}^*$ for some $p \geq 0$.
It is well known \cite{Lag}
that $U_0 = LP$, where $LP$ is the Laguerre-P\'olya class of
entire functions which are locally uniform limits of real polynomials
with real zeros.

The following results established conjectures of Wiman \cite{Al1,Al2}
and P\'olya \cite{Polya43} respectively.
Here all counts of zeros should be understood to
be with respect to multiplicity, and the same convention will be
maintained throughout the paper unless explicitly stated otherwise.

\begin{thm}[\cite{EdwH,SS}]\label{thmC}
Let $p \in \N$ and let $f \in U_{2p}^*$. Then $f''$
has at least $2p$ non-real zeros.
\end{thm}

\begin{thm}[\cite{BEpolya}]\label{thmB}
Let $p$ be a positive integer and let $f \in U_{2p}^*$. Then the number
of non-real zeros of the $k$th derivative
$f^{(k)}$ tends to infinity with $k$.
\end{thm}
\begin{thm}[\cite{BEL,lajda}]\label{thmA}
If $f$ is
a real entire function of infinite order 
then $ff^{(k)}$ has infinitely many
non-real zeros, for every $k \geq 2$.
\end{thm}

The present paper addresses the following related problem: if 
$f$ is a real entire function with finitely many non-real zeros,
must a linear differential polynomial
$$\Psi = f^{(k)} + a_{k-1} f^{(k-1)} + \ldots + a_0 f$$
with
constant real coefficients $a_j$ have non-real zeros, and if so
how many? This question will be
resolved for $k=2$, in which case
in view of the standard transformation 
$$f(z) = e^{-a_1 z/2} g(z), \quad
\Psi (z) = e^{-a_1 z/2} ( g''(z) + (a_0 - a_1^2 /4) g ),$$ 
it may be assumed with no loss of generality that $a_1 = 0$.

\begin{thm}\label{thm1}
Let $f$ be a real entire function 
with  finitely many non-real zeros, and let 
$\omega $ be a positive real number. If
$f \in U_{2p}^*$
for some $p \in \N$ then $f'' + \omega f$ has at least $2p$ non-real zeros.
If $f$ has infinite order then
$f'' + \omega f$ has infinitely many non-real zeros.
\end{thm}


It evidently suffices to prove Theorem \ref{thm1}
for $\omega = 1$, but the following examples
show that the theorem fails for $\omega < 0$. If $f$ is defined by
\cite{FH}
$$
\frac{f'(z)}{f(z)} = a + e^{-2az}, \quad
\frac{f''(z)}{f(z)} = a^2 + e^{-4az} ,
$$
then $f$ and $f'' - a^2 f$ have no zeros
at all in the plane. For an example of finite order define a 
zero-free function
$f \in U_2$ by setting \cite{BEL2}
$$
\frac{f'(z)}{f(z)} = -16 z^2 + 8z + 2, \quad
\frac{f''(z) - 12 f(z)}{f(z)} = 256 z^3 (z-1),
$$
so that $f''-12f$ has only real zeros.

The proof of Theorem \ref{thm1} will use machinery developed in
\cite{LeO,SS} for the Wiman conjecture, and refinements from
\cite{BEpolya,BEL,lajda}, but will depart from the earlier methods
in several significant steps. The aim is to construct an auxiliary 
function having finitely many critical points in $\C \setminus \R$, and 
this will be done in Lemma \ref{lem1}, but in contrast to
\cite{BEpolya,BEL,lajda,LeO,SS} the resulting function may have
a finite non-real asymptotic value. Moreover
for the present problem the normal families
arguments used successfully in \cite{BEL,lajda} seem difficult to
apply, since the condition
$$
f(z) (f''(z) + f(z)) \neq 0 
$$
is not invariant under a change of variables $w = Rz$. It also seems
worth observing that for $f$ in $U_{2p}^*$,
whereas every derivative of $f$ has finitely many  non-real zeros
(see e.g. \cite[Corollary 2.12]{EdwH}), this need not be the case for
$f'' + f $, as the simple example
$f(z) = 1 + \sin (z/2) \in U_0$ shows.
For further remarks and contrasts see \S\ref{contrast}.

\section{Preliminaries}
\begin{definition}\label{def1}
For $a \in \C$ and $0 \leq s < r < R \leq + \infty $ set
$$
D(a, r) = \{ z \in \C : | z- a | < r \}, \quad
S(a, r) = 
\partial D(a, r),
\quad A(s, R ) = \{ z \in \C 
: s < |z| < R \} 
$$
and 
$$
H = \{ z \in \C : {\rm Im} \, z > 0 \}, \quad
D^+(0, r) = D(0, r) \cap H, \quad
A^+(s, R ) = A(s, R ) \cap H .
$$
\end{definition}

\begin{lem}[\cite{Tsuji}]\label{lemhm}
Let $u$ be a non-constant continuous subharmonic function in the plane.
For $r > 0$
let $\theta^*(r)$ be the angular measure of that subset of $S(0, r)$
on which $u(z) > 0$, except that $\theta^*(r) = \infty $ if 
$u(z) > 0$ on the whole circle $S(0, r)$. Then, for $r > 0$,
\begin{equation*}
B(r, u) = \max \{ u(z) : |z| = r \}
\leq 3 T(2r, u) = \frac3{2 \pi} 
\int_0^{2 \pi} \max \{ u( 2r e^{it} ) , 0 \} \, dt 
\end{equation*}
and, if $r \leq R/4$ and $r$ is sufficiently large, 
\begin{equation*}
B(r, u) \leq 9 \sqrt{2} B(R, u) 
\exp \left( - \pi \int_{2r}^{R/2} \frac{ds}{s \theta^*(s)} \right).
\end{equation*}
\end{lem}
\hfill$\Box$
\vspace{.1in}

\begin{lem}\label{lemE1}
Let $0 < A < B < \infty$ and $0 < M < \infty$, and suppose that
$D_1 , D_2, \ldots, D_N$ are pairwise disjoint simply connected
domains, each lying in $\C \setminus \{ 0 \}$ and satisfying
$$
\int_{r^A}^{r^B} \frac{ \pi \, dt}{t \theta_{D_j}(t) } \leq M \log r ,
$$
where $\theta_{D_j}(t)$ denotes the angular measure of $D_j \cap
S(0, t)$. Then $N (B-A)\leq 2M$.
\end{lem}
\textit{Proof.} This is completely
standard. The Cauchy-Schwarz inequality gives
$$
N^2 = \left( \sum_{j=1}^N 1 \right)^2 
\leq \left(
\sum_{j=1}^N \theta_{D_j}(t) \right) \left( \sum_{j=1}^N \frac1{\theta_{D_j}(t)}
\right)
\leq 2 \pi \sum_{j=1}^N \frac1{\theta_{D_j}(t)} . 
$$
Hence
$$
N^2 (B-A) \log r \leq   \sum_{j=1}^N
\int_{r^A}^{r^B} \frac{ 2 \pi \, dt}{t \theta_{D_j}(t) } \leq  2N M \log r .
$$
\hfill$\Box$
\vspace{.1in}

The proof of Theorem \ref{thm1} requires
the characteristic function in a half-plane as developed in
\cite{LeO,Tsuji0} 
(see also \cite{BEL,GO}).
Let $g$ be meromorphic
in a domain containing the closed upper half-plane $\overline{H} =
\{ z \in \C : \mathrm{Im} \, z \geq  0 \}$.
For $t \geq 1$ let
$\mathfrak{n} (t, g) $ be the number of poles of $g$
in $\{ z:|z-it/2|\leq t/2, |z| \geq 1\}$, and 
for $ r \geq 1$ set
\begin{equation}
\mathfrak{N} (r, g) = \int_1^r
\frac{ \mathfrak{n} (t, g) }{t^2} \, dt, \quad
\mathfrak{m} (r, g) =
\frac1{2 \pi} \int_{ \sin^{-1} (1/r)}^{ \pi - \sin^{-1} (1/r)}
\frac{ \log^+ | g(r \sin \theta e^{i \theta } )|}
{r \sin^2 \theta } \,  d \theta.
\label{tsujim}
\end{equation}
The Tsuji characteristic $\mathfrak{T} (r, g)$ is then given by
$ \mathfrak{T} (r, g)  = \mathfrak{m} (r, g)  +\mathfrak{N} (r, g)$.
\begin{lem}[\cite{LeO}]\label{lemLO}
Let $g$ be meromorphic in $\overline{H}$ such that 
\begin{equation*}
\mathfrak{m}(r, g) = O( \log r ) \quad\hbox{as}\quad r \to \infty ,
\end{equation*}
where $\mathfrak{m}(r, g) $ is given by $(\ref{tsujim})$.
Then, as $R \to \infty$,
$$
\int_R^\infty \frac{m_{0\pi} (r, g) }{r^3} \, dr 
\leq
\int_R^\infty \frac{\mathfrak{m} (r, g ) }{r^2} \, dr
= O \left( \frac{ \log R}{R} \right),  \quad
m_{0\pi} (r, g) =
\frac1{2\pi} \int_0^\pi \log^+ |g(r e^{i \theta  } )| \, d \theta. 
$$
\end{lem}
\hfill$\Box$
\vspace{.1in}

The next lemma involves direct transcendental singularities of the
inverse function \cite{BE,Nev}.
Let $a \in \C$ be an asymptotic value of the transcendental meromorphic 
function $g$, so that $g(z) \to a$ as $z \to \infty$ along a path
$\gamma $ tending to infinity. Then the inverse function $g^{-1}$ is
said to have a transcendental singularity
over $a$. For each $\varepsilon > 0$ 
there exists a component $C = C( a, \varepsilon, g)$ of the
set $\{ z \in \C : |g(z) - a | < \varepsilon \}$ with the property that
$C$ contains an unbounded subpath of $\gamma$. Two asymptotic paths 
$\gamma, \gamma'$ on which $g(z) \to a$ 
determine distinct singularities if the corresponding components
$C( a, \varepsilon, g)$, $C'( a, \varepsilon, g)$ are distinct for  
some $ \varepsilon > 0$.

The singularity of $g^{-1}$ corresponding to $\gamma$
is called indirect if $C( a, \varepsilon, g)$, for every
$\varepsilon > 0$, contains infinitely many zeros of $g-a$ \cite{BE}, and 
direct otherwise, in which case
$C( a, \varepsilon, g)$, for all sufficiently small
$\varepsilon > 0$, contains no zeros of $g-a$. 
With a slight abuse of
notation, such a singularity will be referred to as lying in the
upper half-plane $H$ if
$C(a, \varepsilon, g) \subseteq H$ for sufficiently small positive
$\varepsilon $. Transcendental singularities over $\infty$ are defined and
classified analogously. 

\begin{lem}\label{directsinglem}
Let $g$ be a meromorphic function in the plane such that
$\mathfrak{T}(r, g) = O( \log r )$ as $r \to \infty$. Then there is 
at most one direct singularity of $g^{-1}$ lying in $H$.
\end{lem}
\textit{Proof.} Assume that $g^{-1}$ has at least two
direct singularities over $a_1, a_2$ in $H$. Here $a_1, a_2$ need not be
distinct but may be assumed finite. Hence for some $\varepsilon  \in (0, 1)$
and for $j=1, 2$ 
there exists an unbounded component $D_j \subseteq H$ of the set 
$\{ z \in \C  : |g(z) - a_j| < \varepsilon \} $, such that $g(z) \neq a_j$
on $D_j$. The functions $u_1, u_2$ defined by
$$
u_j(z) = \log |\varepsilon /(g(z) - a_j) | \quad (z \in D_j), \quad
u_j(z) = 0 \quad (z \not \in D_j),
$$
are then non-constant and subharmonic in the plane with disjoint
supports $D_j \subseteq H$. Since $\mathfrak{T}(r, 1/(g - a_j)) =
O( \log r )$ as $r \to \infty$, Lemmas \ref{lemhm} and
\ref{lemLO} lead to, for large positive $R$,
\begin{equation}
\frac{B(R/2, u_j)}{2 R^2}
\leq \int_R^\infty \frac{B(r/2, u_j)}{r^3} \, dr \leq
3 \int_R^\infty \frac{m_{0\pi}(r, 1/(g-a_j))}{r^3} \, dr = O
\left( \frac{ \log R}{R} \right) ,
\label{dir0}
\end{equation}
and hence $B(R, u_j) = O(R \log R)$ as $R \to \infty$. But applying Lemma 
\ref{lemhm} again and using the Cauchy-Schwarz inequality as
in the proof of Lemma \ref{lemE1}, as well as the fact that
the $D_j$ are disjoint and lie in $H$, now yields
$$
4 \leq \pi \sum_{j=1}^2 \frac1{ \theta_{D_j} (s) }, \quad
4 \log R \leq
\int_1^R  \sum_{j=1}^2 \frac{ \pi \, ds}{s \theta_{D_j} (s) } \leq
(2 + o(1) ) \log R 
$$
as $R \to \infty$, which is plainly a contradiction.
\hfill$\Box$
\vspace{.1in}

The following lemma is the well known Carath\'eodory inequality
\cite[Ch. I.6, Theorem $8'$]{Le} for analytic self-mappings of 
the upper half-plane $H$.

\begin{lem}
Let $\psi : H \to H$ be analytic. Then
\begin{equation}
\frac{|\psi(i)|\sin\theta}{5r}<|\psi(re^{i\theta} )|<
\frac{5 r |\psi(i)|}{\sin\theta} \quad \hbox{for}  \quad r \geq 1,\,
\theta\in (0,\pi).
\label{cara}
\end{equation}
\end{lem}
\hfill$\Box$
\vspace{.1in}

The proof of Theorem \ref{thm1} will require some elementary
inequalities for the hyperbolic metric on
the upper half-plane $H = \{ z = x + iy: x \in \R, y > 0 \}$, on which 
the hyperbolic density is $1/y$. Hence if $\gamma$ is a curve
joining $i$ to $z = x + iy \in H$ then the hyperbolic length of 
$\gamma$ is
\begin{equation}
[i, z]_H  =
\int_\gamma \frac1{{\rm Im} \, \zeta } \, | d \zeta |
  \geq \left| \int_1^y \frac1t \, dt \right| 
= | \log y | = \left| \log \left( \frac1{{\rm Im}\, z }  \right) \right| .
\label{hyp1}\end{equation}
On the other hand $i$ may be joined to $z$ by the line segment
$\gamma_1$ from $i$ to $x + i$ followed by the line segment
$\gamma_2$ from $x+i$ to $z$, which gives the upper bound
\begin{equation}
[i, z]_H \leq \left( \int_{\gamma_1} + \int_{\gamma_2} \right) 
\frac1{{\rm Im} \, \zeta } \, | d \zeta |
\leq | {\rm Re} \, z | + \left| \log \left( \frac1{{\rm Im}\, z }  \right)
\right| .
\label{hyp2}\end{equation}
The imaginary parts of $T = \tan z$ and $z$ will now be compared for
$z \in H$. It is clear that
$T = \tan z = M(u)$
for $z \in H$, where $M : D(0, 1) \to H$ is a M\"obius transformation
with $M(0) = i$,
and $u =  e^{2iz}$ maps $H$ into $D(0, 1)$. If 
${\rm Im}\, T$ is small then evidently so is $y = {\rm Im}\, z$, and
$$
\log \left( \frac{1+|u|}{1-|u|} \right) =
[0, u]_{D(0,1)} = [i, T]_H \geq 
\log \left( \frac1{{\rm Im}\, T }  \right),
$$
using (\ref{hyp1}). Hence
\begin{equation}
2y \sim 1- e^{-2y} = 1- |u| \leq  2 \, {\rm Im}\, T, \quad T =  \tan z ,
\label{hyp3}\end{equation}
uniformly in $x = {\rm Re}\, z $ as $y = {\rm Im}\, z$ tends to $0$.
\hfill$\Box$
\vspace{.1in}

\section{Direct singularities and critical points}\label{rh}

If an analytic function is a proper mapping between 
domains each of finite connectivity then the Riemann-Hurwitz formula 
\cite[p.7]{Stei2} links the valency of the mapping with the number of
critical points and the connectivities of the domains. To apply this
formula requires that the function map boundary to boundary in the sense
of \cite[p.4]{Stei2}.

The function $f(z) = z e^z $ has a direct transcendental
singularity over $0$, and a critical point at $-1$, 
and the interval $(- \infty, 0]$ lies in a component $C$ of the set
$\{ z \in \C : |f(z)| < 1 \}$. The function $f$ is infinite-valent
on $C$, but the number of zeros of $f$ in $C$ is equal
to the number of critical points of $f$ in $C$.
The proof of Theorem \ref{thm1} will require a relation between
zeros and critical points for components of this type,
and this will be obtained by transforming the function to one of form
$R(z) \exp( az )$ with $R$ a rational function and $a \in \C$.

\begin{lem}\label{explem}
Let $b $ be a positive real number and let $R$ be a rational
function such that $|R(x)| = 1$ for all $x \in \R$. Assume that
$f(z) = R(z) e^{i bz} $ is such that $f$ has no critical values $w$
with $|w| = 1$. Let $A \subseteq H$ be an unbounded component of the set
$\{ z \in \C : |f(z)| < 1 \}$,
and let $p$ be the connectivity of $A$. Let $m$ be the number of
zeros of $f$ in $A$ and $n$ the number of zeros of $f'$ in $A$. Then
$m - n = 1-p$.
\end{lem}
\textit{Proof.}
It is evident that such a component $A$ exists, because
$|f(iy)| < 1$ for all large positive real $y$ and
$|f(x)| = 1$ for $x \in \R$.
The set $X = \{ z \in \C : |f(z)| = 1\}$ consists of pairwise disjoint 
Jordan curves, and
Jordan arcs tending to infinity in both directions, one of
which is the real axis. Since, as $z = r e^{i \theta} \to \infty$, 
$$
\log |f(re^{i \theta} )| = - br \sin \theta + O(1), \quad
\frac{ \partial \log |f(re^{i \theta} )| }{ \partial \theta } =
- {\rm Im} \, \left( \frac{zf'(z)}{f(z)} \right) =
- br \cos \theta + o(1),
$$
it follows that if $z \in X$ is large then $z \in \R$. Hence
the finite boundary $\partial A$ consists of the real axis and $p-1$
pairwise disjoint Jordan curves $\Gamma_j$ in $H$. Let
$\overline{\Gamma}_j$ be the reflection of $\Gamma_j$ in the real axis.
Let $t$ be large and positive and let $\gamma$ be the cycle
consisting of the circle $S(0, t)$ described once counter-clockwise
and each of the $\Gamma_j$ and $\overline{\Gamma}_j$ described
once clockwise. Since $f(z) = 0$ if and only if $f ( \overline{z}) =
\infty$, and since the multiplicities coincide, the net change
in $\arg f (z)$ as $z$ describes $\gamma$ is $0$.

Because $t$ is large it follows that $f'(z) \sim ib f(z)$ on
$S(0, t)$ and so the net change in 
$\arg f' (z)$ as $z$ describes $S(0, t)$ agrees with that of
$\arg f(z)$. Moreover, the net change in 
$\arg f' (z)$ as $z$ describes one of the 
$\Gamma_j$ or $\overline{\Gamma}_j$ clockwise exceeds that of 
$\arg f(z)$ by $2 \pi $ \cite[p.122]{Titch}. Hence
if $N$ is the number of zeros minus the number of poles of $f'$
which lie inside $\gamma$ (i.e. which have
winding number $1$ relative to $\gamma$),
then $N = 2 (p-1)$.

Now the only zeros and poles of $f$ which lie inside $\gamma$ are the
the zeros of $f$ in $A$ and their reflections across
$\R$, which are poles. Let these zeros of $f$ in $A$ be denoted by
$z_j$, with multiplicities $p_j$. Then 
$z_j$ and $\overline{z}_j$ 
together contribute $p_j - 1 - (p_j + 1) = -2$ to $N$.
Next let $w_k$ be the zeros of $f'$ in $A$ which are not zeros of $f$, 
and denote their multiplicities by $q_k$. Then 
$w_k$ and $ \overline{w}_k$ together contribute 
$2 q_k$ to $N$. 
Let $r$ be the number of distinct zeros $z_j$ of $f$ in $A$. Then
summing over the $z_j$ and $w_k$ gives
$$
2(m-n) = 2 \left( r - \sum q_k \right) = -N = 2(1-p) .
$$
\hfill$\Box$
\vspace{.1in}

Recall next some
standard facts from \cite[p.287]{Nev}, albeit in
slightly more general form. 
Let the function $G$ be transcendental and meromorphic in the plane,
with no asymptotic values in
$$
V_1 = \{ v \in \C : 0 < |v| < 1 \} ,
$$
and assume further that $G'$ has finitely many zeros $z$ with 
$G(z) \in V_1$.  
Let $\Gamma$ be a simple piecewise analytic arc, starting at $v_1 
\in S(0, 1)$
but otherwise lying in $V_1$, such that all critical values 
$v \in V_1$ of $G$ lie on $\Gamma$. Choose a branch of the
logarithm defined near to $v_1$ and let $\gamma = \log \Gamma $,
so that $\gamma$ is a simple piecewise analytic arc
and $e^\gamma = \Gamma$. For $k \in \Z$ let
$\gamma_k $ be the translation by $k 2 \pi i$ of $\gamma$; these
$\gamma_k$ are then pairwise disjoint. 
Now let 
\begin{equation}
V_0 = V_1 \setminus \Gamma, \quad
U_0 = K(0) \setminus \bigcup_{k \in \Z} \gamma_k ,
\quad \hbox{where} \quad 
K(t) = \{ u \in \C : {\rm Re} \, u < t \} .
\label{rh1}
\end{equation}
Then $\exp (U_0) =V_0$. Let $C$ be a component
of $G^{-1}(V_0)$, and choose $z_0 \in C$ and $u_0 \in U_0$ with
$G(z_0) = v_0 = \exp( u_0)$. Let $g$ be the branch of $G^{-1}$ 
mapping $v_0$ to $z_0$. Then
\begin{equation}
h(u) =  g( e^u) =  G^{-1}( e^u) 
\label{rh2}
\end{equation}
extends by the monodromy theorem to be analytic on $U_0$, with
$h(U_0) \subseteq C$. Indeed if $z \in C$ then $z_0$ may be joined
to $z$ by a path $\lambda_1$ in $C$ and there exists a path
$\lambda_2$ in $U_0$ starting at $u_0$ such that
$\exp( \lambda_2) = G( \lambda_1) \subseteq V_0$. 
Then $\lambda_1 = h ( \lambda_2)
\subseteq h(U_0)$, since $\lambda_1 $ and $ h ( \lambda_2)$ both
start at $z_0$ and have the same image under $G$. Hence
$h(U_0) = C$.

Suppose first that $h$ is univalent on $U_0$. Then for $t < 0$ with
$|t|$ large the image of
the line ${\rm Re} \, u = t$ under $h$ is a level curve 
$|G(z)| = e^t$ which tends to infinity in both directions. Hence
$h(u) \to \infty $ as $u \to \infty $ in $K(t)$, and $C$ is an unbounded
simply connected domain containing a path tending to infinity on
which $G(z) \to 0$. Such components of $G^{-1}(V_0)$ will be called
type I.

If the finite boundary $\partial C$ of a type I component
contains no critical point $z$ of $G$ with $|G(z)| < 1$
then $h$ may be continued analytically
along each $\gamma_k$ to be univalent on $K(0)$,
and $C$ lies in a component 
$B = h(K(0))$ of $\{ z \in \C : |G(z)| < 1 \}$ which contains no zeros
of $G$.

Suppose next that $h$ is not univalent on $U_0$.
Then there exist distinct $u_1, u_2 \in U_0$ 
with $h(u_1) = h(u_2)$ and hence $e^{u_1} = e^{u_2}$. Take the least $k \in \N$ 
for which there exist $u_3, u_4 \in U_0$ 
with $u_3 = u_4 + k 2 \pi i $ and $h(u_3) = h(u_4)$. Then
$h$ has period $k 2 \pi i$ by
the open mapping theorem and
$$
F( \zeta ) = g( \zeta^k ) = h( k \log \zeta ) 
$$
extends to be analytic in 
$Z_k = \{ \zeta \in \C : \zeta^k \in V_0 \},$
mapping $Z_k$ univalently onto $C$. Moreover, 
$z_1 = \lim_{\zeta\to 0} F(\zeta)$
exists, and must be finite, since otherwise every large
$z \in \C$ is $F(\zeta)$ for some $\zeta \in Z_k$ and satisfies
$G(z) = \zeta^k \in V_0$, contradicting the assumption that
$G$ is transcendental. Hence $z_1$ is a zero of $G$ 
and $G$ maps $C \cup \{ z_1 \}$ onto
$V_0 \cup \{ 0 \}$, the mapping $k$-valent. 

This time $C$ will be called type II. Here
if $\partial C$ contains no zero $z$ of $G'$ with
$0 < |G(z)| < 1$ then
$F$ may be analytically continued to $D(0, 1)$, with the
extended function univalent by the open mapping theorem, and $C$
lies in a component $B = F(D(0, 1))$ of $\{ z \in \C : |G(z)| < 1 \}$ 
which contains the zero $z_1$ and is such that $G$ is $k$-valent
on $B$.

Now let $A$ be any component of the set $\{ z \in \C : |G(z)| < 1 \}$
and let $C \subseteq A$ be a component of the set $G^{-1}(V_0)$. If
$\partial C$ contains no zero $z$ of $G'$ with $0 < |G(z)| < 1$ then 
$C$ is the only component of $G^{-1}(V_0)$ contained in $A$, and
$G$ has at most one zero in $A$, possibly multiple. In the general case,
it follows from the fact that $G'$ has finitely
many zeros $z$ with $0 < |G(z)| < 1 $ that $A$ contains
finitely many components $C$ of $G^{-1}(V_0)$ and finitely
many zeros of $G$. Moreover if $A$ does not contain any type I 
components $C$ of $G^{-1}(V_0)$ nor any zeros of $G'$ then 
$A$ contains one simple zero of $G$ and $G$ is univalent on $A$.

\begin{lem}\label{rhlem}
With $G$ as above and $V_0$ defined as in $(\ref{rh1})$
let $A$ be a component of the set $\{ z \in \C : |G(z)| < 1 \}$ 
containing precisely one
type I component $C$ of $G^{-1}(V_0)$. Then the number of
zeros of $G$ in $A$ is at most the number of zeros of $G'$ in $A$.
\end{lem}
\textit{Proof.} Choose $z_0 \in C$ such that $t = |G(z_0)|$ is small, 
and join $z_0$ to each zero of $G$ in $A$ 
by a path in $A$. The union of
these finitely many paths forms a compact connected set $E \subseteq A$
with
$$
\max \{ |G(z)| : z \in E \} < 1 ,
$$
and $E$ is contained in a component $\widetilde A \subseteq A$ of the set 
$\{ z \in \C : |G(z)| < 1 - \delta  \}$ for some small positive
$\delta $. Set $\widetilde G = G/(1-\delta)$. Then a set
$\widetilde V_0$ may be defined corresponding to $\widetilde G$ in the same way
as $V_0$ was defined for $G$, and 
since $C$ contains a path tending to infinity
on which $G(z) \to 0$ it is clear that
$\widetilde A$ contains at least one type I component of
$\widetilde G^{-1} (\widetilde V_0 )$. 

Suppose on the other hand that 
$W_1, W_2$ are distinct type I components of
$\widetilde G^{-1} (\widetilde V_0 )$ contained in $\widetilde A$. 
Choose $w_j \in W_j$ with $G(w_j)$ small and hence $w_j$ large. Then
$w_1, w_2$ must both lie in $C$ and may be joined in $C$ by a path 
$\sigma$ on which $G(z)$ is small. But then $\sigma $ lies in a 
component of $\widetilde G^{-1} (\widetilde V_0 )$ and this is a contradiction.

These observations show that
in order to prove Lemma \ref{rhlem}
there is no loss of generality in assuming
there exists a small positive $\eta$ such
that $G$ has no asymptotic values $v$ with $0 < |v| \leq 1 + \eta$, 
and that $G'$ has finitely many zeros $z$ with $0 < |G(z)| \leq 1 + \eta$, and
none with $|G(z)| = 1$, since
otherwise $G$ may be replaced by~$\tilde G$.

Choose $u_0 \in U_0$ with $\exp( u_0) = v_0 = G(z_0)$ and define $h$ as in
(\ref{rh2}) using the branch of $G^{-1}$ mapping 
$v_0$ to $z_0$. Since $h$
extends to be univalent on $U_0$ and $G'$ has finitely many zeros $z$ with
$0 < |G(z)| \leq 1$, it follows that if $|k|$ is large then $h$ may be 
continued along the arc $\gamma_k$ and the extended function is
still univalent. Indeed, if $S$ is large enough then $h$ extends to 
be analytic and univalent on the set
$U_1 = \{ u \in \C : {\rm Re} \, u \leq 0, |u| \geq S \}$.

Since there are no asymptotic values $v$ of $G$ with $0 < |v| \leq 1$,
all type II components of $G^{-1}(V_0)$ are bounded, and since there
are finitely many of these contained in $A$, say
$D_j$, it follows that there exists
$R > 0$ such that $E$ and all the $D_j$ lie in $D(0, R)$. Moreover
$R$ may be chosen so large that $|h(u)| < R$ for all 
$u \in U_0 \cap D(0, 2S)$.

The components of the finite boundary $\partial A$ 
are pairwise disjoint level curves $|G(z)| = 1$,
each either a Jordan curve or a Jordan arc tending to infinity
in both directions.
If $\Lambda$ is an unbounded component of 
$\partial A$ then each large $z \in \Lambda $ must belong to
$\partial C$ and so must be $h(is)$ for some real $s$ with $|s|$ large.
Hence there is precisely one unbounded component $\Lambda$ of
$\partial A$. 
Moreover all but finitely many $1$-points of $G$ in $\partial A$
lie on $\Lambda$ and $\partial A$ has finitely many components.

Let $\Omega$ be the component of $\C \setminus \Lambda$ which contains
$A$, and let $z = p(w)$ map the upper half-plane $H$ conformally onto
$\Omega$. Then the function $q$ defined by
$$
q(w) = G(p(w)) \quad (w \in H), \quad \overline{q(\overline{w})} =
\frac1{q(w)} ,
$$
extends by the reflection principle to a meromorphic function on the plane,
which must have the form $q(w) = R(w)e^{i S(w)} $ with $R$ a rational
function such that $R(\infty) = 1$, and $S$ an entire function which must
be real since $|q(w)| =
|R(w)| = 1$ on $\R$. Moreover, if $w$ is large and $|q(w)| = 1$, 
then $w $ is real, since $\partial A$ has finitely many components, of
which only $\Lambda$ is unbounded. Since $A$ contains a path 
tending to infinity on which $G(z) \to 0$, it follows 
that $|q(w)| < 1$ and ${\rm Re} \, ( i S(w) ) \leq o(1)$
for all large $w \in H$.
It now follows from the Wiman-Valiron theory \cite{Hay5}
that $S$ is a polynomial, which must be of form
$S(w) = aw+b$ with real
constants $a, b$ and $a > 0$. 
Since $p^{-1}(A) \subseteq H$ is a component of the set
$\{ w \in \C : |q(w)| < 1 \}$, 
the result now follows from Lemma \ref{explem}.
\hfill$\Box$
\vspace{.1in}

\section{Proof of Theorem \ref{thm1}: first steps}

Let $f$ be a real transcendental
entire function and assume that
$f$ and $f'' +  f$ have finitely many non-real zeros, and that
either $f$ has infinite order 
or $f \in U_{2p}^*$ for some positive
integer $p$. 

\begin{lem}\label{lemA1}
Set $L = f'/f$. Then $L$ satisfies 
\begin{equation}
\mathfrak{T}(r, L) = O( \log r ) \quad \hbox{as $r \to \infty$}.
\label{A2}
\end{equation}
\end{lem}
\textit{Proof.} This uses a modified Tumura-Clunie argument
\cite[p.69]{Hay2} (see also \cite{HellY}). Write
\begin{equation}
M = L' + L^2 + 1   = \frac{f''+f}{f} , \quad M' = L'' + 2 L L' = 
\frac{M'}{M} ( L' + L^2 + 1).
\label{A3}
\end{equation}
Then
\begin{equation}
2PL = Q = \frac{M'}{M} ( L' + 1) - L'' ,
\quad \hbox{where} \quad P = L' - \frac{M'}{2M} L .
\label{A4}
\end{equation}
But $L$ has finitely many non-real poles, and $M$ has finitely many
non-real zeros. Since the lemma of the logarithmic derivative holds for
the Tsuji characteristic \cite[p.332]{LeO}, so does
a direct analogue
of Clunie's lemma \cite[p.68]{Hay2}, which on combination with
(\ref{A3}) gives
\begin{equation}
\mathfrak{T}(r, P) + \mathfrak{T}(r, M'/M) = \mathfrak{S}(r, L),
\label{A5}
\end{equation}
where $\mathfrak{S}(r, L)$ denotes any quantity which satisfies
$$
\mathfrak{S}(r, L) \leq o( \mathfrak{T}(r, L) ) + O( \log r )
$$
as $r \to \infty$, possibly outside a set of finite measure. Now write
\begin{equation}
U = L + \frac{M'}{4M} , \quad
M = L^2 + \frac{M'}{2M} L + P + 1 = U^2 + R, \quad
\mathfrak{T}(r, R) = \mathfrak{S}(r, L),
\label{A6}
\end{equation}
using (\ref{A3}), (\ref{A4}) and (\ref{A5}). Thus
\begin{equation}
M' = 2UU' + R' , \quad
UV = \frac{M'}{M} R - R', \quad
\hbox{where} \quad V = 2U' - \frac{M'}{M} U ,
\label{A7}
\end{equation}
and, using (\ref{A6}) and (\ref{A7}) and Clunie's lemma,
\begin{equation}
\mathfrak{T}(r, U) = \mathfrak{T}(r, L) + \mathfrak{S}(r, L), \quad
\mathfrak{T}(r, V) +  \mathfrak{T}(r, UV) = \mathfrak{S}(r, L) .
\label{A8}
\end{equation}
If $V \not \equiv 0$ then (\ref{A8}) gives
$$
\mathfrak{T}(r, L) \leq \mathfrak{T}(r, U) + \mathfrak{S}(r, L)
\leq \mathfrak{T}(r, UV) + \mathfrak{T}(r, V) +  \mathfrak{S}(r, L)
=  \mathfrak{S}(r, L),
$$
which gives (\ref{A2}). Assume henceforth that $V \equiv 0$. Then
there exists a constant $d$ such that
\begin{equation}
M = dU^2 , \quad (d-1)U^2 = R ,
\label{A9}
\end{equation}
and it may be assumed that $d = 1$, since otherwise
(\ref{A6}) and (\ref{A8}) give (\ref{A2}). Thus, by (\ref{A6}),
\begin{equation}
L = W + cM^{1/2}, \quad \hbox{where} \quad W =
-  \frac{M'}{4M} \quad \hbox{and} \quad c^2 = 1.
\label{A10}
\end{equation}
This gives
$$ 
L' = W' + 
\frac12 c M^{-1/2} M' = W' - 2WcM^{1/2}
$$
and 
$$
M = L^2 + L' + 1 = 
M + 2WcM^{1/2} + W^2 + L' + 1 = M + W^2 + W' + 1.
$$
It follows using (\ref{A9}) and (\ref{A10}) that
$$
0 = W^2 + W' + 1, \quad
W(z) = - \tan (z + A), \quad 
U(z) = B \sec^2 (z + A), \quad A, B \in \C.
$$
But $e^{iz}$ is bounded in $H$ and so 
$\mathfrak{T}(r, U) = O(1)$, from which (\ref{A2}) follows
using (\ref{A8}) again.
\hfill$\Box$
\vspace{.1in}

The next step is the Levin-Ostrovskii factorisation \cite{LeO} 
of $L = f'/f$, which will be developed following \cite{SS} but using
refinements from \cite{BEpolya}, slightly modified.

\begin{lem}\label{factorlem}
The logarithmic derivative $L = f'/f$ has a factorisation
\begin{equation}
L = \frac{f'}{f} = \phi \psi  
\label{B1}
\end{equation}
in which $\phi$ and $\psi$ are real meromorphic functions
satisfying the following:\\
(i) 
either $\psi \equiv 1$ or $\psi (H) \subseteq H$;\\
(ii) $\psi$ has a simple pole at each real zero of $f$, and no other poles;\\
(iii) $\phi$ has  finitely many poles, none of them real;\\
(iv) on each component of $\R \setminus f^{-1}( \{ 0 \} )$ the number of
zeros of $\phi$
is either infinite or even;\\
(v) if $f \in U_{2p}^*$ then $\phi$ is a rational function, and
if in addition $f$ has at
least one real zero then 
the degree at infinity of $\phi$ is even and satisfies
\begin{equation}
{\rm deg}_\infty (\phi) = 
\lim_{z \to \infty} \frac{ \log | \phi(z) | }{\log |z| }  \geq 2p.
\label{degphidef}
\end{equation}
\end{lem}
Here a meromorphic function $g$ on $\C$ is called real if
$g(\R) \subseteq \R \cup \{ \infty \}$.\\
\textit{Proof.} 
Suppose first that $f$ has no real zeros. Then $L = f'/f$
has finitely many poles, and if $f \in U_{2p}^*$ then $L$ is a 
rational function by the lemma of the logarithmic derivative. 
If the number of real zeros of $L$ is infinite or even, set
$\psi = 1$ and $\phi = L$. On the other
hand if $L$ has an odd
number of real zeros $b$, choose such a zero $b$
and write $\psi (z) = z - b$
and $\phi (z) = L(z)/(z - b)$.

Assume henceforth that $f$ has at least one real zero. Then the function
$\psi$ is defined as a product as follows
\cite{BEpolya,SS}. First, 
if $a$ is a real zero of $f$ but not the greatest real zero of $f$,
then there exists a bounded component $(a, b)$ of 
$\R \setminus f^{-1}( \{ 0 \})$. Since $L$ has positive residues
at $a$ and $b$ the number of zeros of $L$ in $(a, b)$ is
odd. Choosing such a
zero $c = c_a \in (a, b)$ of $L$, the factor corresponding to
$a$ is then
\begin{equation}
p_a(z) = \frac{c - z}{a - z} \quad (\hbox{if $ac \leq 0$}), \quad 
\frac{1 - z/c}{1 - z/a} \quad (\hbox{if $ac > 0$}), 
\label{padef}
\end{equation}
and $\arg p_a(z)$ for $z \in H$ is the angle between
the line segments from $z$ to $a$ and $c$ respectively.

Suppose next that $a$ is the greatest real zero of $f$. If
the number of zeros $c$
of $L$ in $(a, \infty)$ is finite but odd, choose such a zero $c$
and form a factor $p_a(z) $
as in (\ref{padef}).
On the other hand if $L$ has an infinite or even number of zeros 
in $(a, \infty)$, take the factor $q_a(z) = 1/(a - z)$, so that
for $z \in H$ the argument $\arg q_a(z)$ is the angle between the 
line segment from $z$ to $a$ and the horizontal line from $z$ in
the direction of $+ \infty$.

Finally, if there is a least real zero $a$ of $f$ and the number of
zeros $c$ of $L$ in $(-\infty, a)$ is finite but odd, then an extra
factor $r_c (z) = z - c$ is included,
and $\arg r_c(z)$ for $z \in H$ is the angle between the line segment from
$z$ to $c$ and the horizontal line from $z$ in
the direction of $- \infty$. 

The function $\psi$ is then the product
of the terms $p_a(z)$ and (if required) $q_a(z)$ and $r_c(z)$, and satisfies
$\arg \psi (z) \in (0, \pi )$ for $z \in H$. Moreover, if
there are infinitely many real zeros $a$ of $f$ then $a c_a > 0$ for
$|a|$ large and so the product converges by the alternating series
test. Furthermore, $\phi$ is defined by (\ref{B1}) and
it is evident from the construction that (i), (ii), (iii)
and (iv) are satisfied.

To establish (v), assume that $f \in U_{2p}^*$ and
recall that by assumption $f$ has at least one real zero
$a_j$. Then (i), (\ref{cara})
and the lemma of the logarithmic derivative give
$$
m(r, \phi) \leq m(r, L) + m(r, 1/\psi) = O( \log r ),
$$
so that $\phi$ is a rational function, using (iii). Since $\phi$ clearly
has an even number of non-real zeros and poles, and has an even number of
real zeros by (iv), the degree at infinity of $\phi$ must be even.
Suppose then that 
\begin{equation}
d_0 = {\rm deg}_\infty  (\phi) \leq 2p - 2, \quad
\phi (z) \sim c_0 z^{d_0} \quad \hbox{as}  \quad z \to \infty, \quad c_0 \neq 0.
\label{frep0}
\end{equation}
Since $\psi (H) \subseteq H$ there
exists $c \geq 0$ such that
\begin{equation}
\psi (z) = c z + o( |z| ) 
\quad \hbox{as}  \quad z \to \infty, \, \pi /4 < \arg z < 3 \pi /4 ,
\label{frep2}
\end{equation}
using the series representation for $\psi$ \cite{BEpolya,Le}.
Combining (\ref{frep0}) and (\ref{frep2}) gives
\begin{equation}
L(z) = c_0 c z^{d_0+1} + o( |z|^{d_0+1} ) \quad \hbox{and} \quad
\log f(z) = \frac{c_0 c z^{d_0+2}}{d_0+2} + o( |z|^{d_0+2} ) =
O( |z|^{2p} )
\label{frep3}
\end{equation}
as $ z \to \infty$ with $\pi /4 < \arg z < 3 \pi /4$.
Write
\begin{equation}
f = P_0 \Pi \exp( P_1 ),
\label{frep1}
\end{equation}
where $P_0$ is a real polynomial with no real zeros, $\Pi$ is
the canonical product formed with the real zeros $a_j$ of $f$,
and $P_1$ is a real polynomial. If $m_j$ is the multiplicity 
of the zero of $f$ at $a_j$ and $A_j $ is the residue of $\psi$
there then, again since $\psi (H) \subseteq H$ \cite{BEpolya,Le},
\begin{equation}
0 < m_j = A_j \phi (a_j), \quad A_j < 0,
\quad \sum_{a_j \neq 0} \frac{ |A_j|}{a_j^2} < \infty .
\label{frep4}
\end{equation}
Hence it follows from (\ref{frep0}) and (\ref{frep4}) that
$$
\sum_{a_j \neq 0} \frac{ m_j }{a_j^{2p}} =
\sum_{a_j \neq 0} \frac{ A_j \phi(a_j) }{a_j^{2p}} \leq
2 |c_0| \sum_{a_j \neq 0} \frac{ |A_j| }{a_j^{2}} + O(1) < \infty .
$$
In particular the product $\Pi$ in (\ref{frep1}) has genus at most
$2p-1$ and growth at most order $2p$, minimal type. 
Since $f \in U_{2p}^*$ the polynomial $P_1$ in (\ref{frep1}) must therefore 
have degree $d_1 \geq 2p$, and if $d_1 = 2p$ then the coefficient
$c_1$ of $z^{d_1}$ in $P_1$ is positive. This gives
\begin{equation}
\log |f(z)| = c_1 {\rm Re} \, (z^{d_1}) + o( |z| ^{d_1} )
\quad \hbox{as}  \quad z \to \infty, \, \pi /4 < \arg z < 3 \pi /4 .
\label{frep5}
\end{equation}
Comparing (\ref{frep3}) and (\ref{frep5}) and recalling that
$d_0 \leq 2p -2$
then forces $d_1 = 2p = d_0 + 2$ and $c_0 c > 0$,
and since $c \geq 0$ both $c_0$ and $c$ must
be positive. Thus $\phi (x) > 0$ for real $x$ with $|x| $ large,
and so $f$ has finitely many real zeros by (\ref{frep4}). 
Hence $L$ and $\psi$ are rational functions and
(\ref{frep2}) holds as
$z \to \infty$ in any fashion.
Moreover, $L(x)$ is positive for large positive $x$, since $ c_0 c > 0$. 
Denoting by $a$ the greatest real zero of $f$ 
it now follows from (\ref{frep0}), (\ref{frep2})
and the fact that $L$ has positive
residue at $a$ that the number of zeros of $L$ in $(a, \infty)$ is
even. 
But then by construction
the product $\psi$ includes a factor $q_a(z) = 1/(a-z)$, so
that $\psi (\infty)$ is finite and $c=0$, contradicting the
conclusion already obtained that $c > 0$. 
\hfill$\Box$
\vspace{.1in}

\begin{lem}\label{lemA2}
The functions $\phi$ and $f$ satisfy
\begin{equation}
T(r, \phi )  +
\log T(r, f) =  O ( r \log r ) ,
\label{phifest}
\end{equation}
as $r \to \infty$, and if $f$ has infinite order
then $\phi$ is transcendental. Moreover, 
there exist $c_1 > 0$ and a set $E_0 \subseteq [1, \infty)$ of
finite logarithmic measure such that
\begin{equation}
\left| \frac{f'(z)}{f(z)} \right| \leq \exp( c_1 r \log r ) 
\label{ggg}
\end{equation}
for large $|z| = r$ outside $E_0$.
Finally, there exists a set $E_1 \subseteq \R$ of measure $0$ 
such that for all $\theta \in \R \setminus E_1$ the estimate
$(\ref{ggg})$ holds as $z = re^{i\theta}$ tends to infinity.
\end{lem}
\textit{Proof.} 
The estimate (\ref{phifest}) and the
fact that $\phi$ is transcendental if $f$ has infinite order
are proved exactly as in \cite[\S6 and \S7]{lajda}
(see also \cite[p.982, pp.989-990]{BEL}):
in particular the bound for $T(r, \phi)$ follows from 
(\ref{cara}), (\ref{A2}), (\ref{B1}) and Lemma \ref{lemLO}. 
The estimates (\ref{ggg}) are now immediate consequences
of standard inequalities
due to Gundersen \cite{Gun2}.
\hfill$\Box$
\vspace{.1in}

\begin{lem}\label{lem1}
Write
\begin{equation}
L = \frac{f'}{f}, \quad T = \tan z , \quad
F = \frac{TL-1}{L+T}.
\label{1}
\end{equation}
Then for any set $X \subseteq \C \setminus \R$
the number of zeros of $F'$ in $X$ is at most the number of
distinct zeros of $f$ in $X$
plus the number of zeros of $f''+f$ in $X$,
and in particular is finite.
Next, let
\begin{equation}
H = \{ z \in \C : {\rm Im} \, z > 0 \}, \quad
W = \{ z \in H : F(z) \in H \}, \quad
Y = \{ z \in H : L(z) \in H \}.
\label{2}
\end{equation}
Then $Y \subseteq W$ and the closure of $Y$ contains no
real zeros of $f$. Moreover if
$C$ is a component of $Y$ then either $\partial C$ contains a non-real zero
of $f$ or $C$ is unbounded and satisfies 
\begin{equation}
\limsup_{z \to \infty, z \in C} {\rm Im} \, L(z) = + \infty .
\label{3}
\end{equation}
Finally,
\begin{equation}
\hbox{$L-i$ and $F-i$ have the same zeros with the same multiplicities.}
\label{L=i}
\end{equation}
\end{lem}
\textit{Proof.} Differentiation of (\ref{1}) gives
$$
F' = \frac{(1+T^2)(L'+L^2 + 1)}{(L+T)^2} =
\frac{(1+T^2)(f'' + f)}{(L+T)^2 f} ,
$$
using (\ref{A3}). Hence non-real zeros of $F'$ can only arise from non-real
zeros of $f''+f$ and non-real zeros of $f$, each of which is a simple pole
of $L$ and hence of $L+T$. It is obvious that a non-real zero of $f''+f$ 
which is not a zero of $f$ is a zero of $F'$ of at most the same
multiplicity. Suppose now that $z$ is a non-real zero of $f$ of 
multiplicity $m \geq 1$, and a zero of $f''+f$ of multiplicity $n \geq 0$.
Since $(L+T)^2$ has a double pole at $z$ it follows that $F'$ cannot 
have a zero at $z$ of multiplicity greater than $n - m + 2$. If $m=1$
this gives $n+1$, which equals the contribution of $z$ to
the number of
distinct non-real zeros of $f$ plus the number of non-real zeros of $f''+f$.
If $m=2$ then $f''(z) \neq 0$ and $n=0$, while if $m \geq 3$ then
$n=m-2$, and both these cases give $n-m+2 = 0$.

Next, (\ref{1}) gives
$$
F = \frac{(TL-1)(\bar L + \bar T)}{|L+T|^2} =
\frac{T |L|^2 + L |T|^2 - \bar L - \bar T}{|L+T|^2} ,
$$
and since $z \in Y$ gives
$$
T \in H, \quad - \bar T \in H, \quad - \bar L(z) \in H,
$$
it follows that $Y \subseteq W$.

Moreover, poles of $L$ coincide with zeros of $f$,
and a real pole of $L$ has
positive residue and so is not in the closure of $Y$ \cite[p.987]{BEL}.
If $C$ is a component of $Y$ such that $\partial C$ does not
contain a non-real zero of $f$ it follows that $C$ is unbounded by the
maximum principle, and the function
\begin{equation}
u_C(z) = {\rm Im} \, L(z) \quad (z \in C), \quad
u_C(z) = 0 \quad (z \not \in  C),
\label{4}
\end{equation}
is non-constant and subharmonic in the plane, from which (\ref{3})
follows.

It remains to prove (\ref{L=i}), which will follow from the fact that
(\ref{1}) gives
$$
F - i = \frac{(L-i)(T-i)}{L+T} .
$$
If $L(z) = i$ then $z$ is non-real and is not a zero of
$L+T$, since $T $ omits the value $-i$, and so $F(z) = i$. Further,
the multiplicities coincide since $T$ omits $i$. Similarly,
if $F(z) = i$ then $z$ is non-real and is not a 
zero of $T-i$ nor a pole of $L$, and so $L(z) = i$.
\hfill$\Box$
\vspace{.1in}

\begin{lem}\label{lemC3}
For $a \in \C \setminus \R$ set
\begin{equation}
s_a = \frac{T(F-a)}{T-F} , \quad T = \tan z .
\label{xx1}
\end{equation}
Then
\begin{equation}
s_a(z) = ( \sin^2 z - a \cos z \sin z  ) L(z) -
\cos z \sin z -  a \sin^2 z 
\label{xx2}
\end{equation}
and $s_a$ has finitely many poles in $H$.

Next, let $M$ and $ N$ be positive real numbers,
and let $b \in \C \setminus \R$ satisfy $b \neq a$. If $z$ is large
with
\begin{equation}
| T - a | > |z|^{-M} \quad \hbox{and} \quad
|F(z) - a | < |z|^{-M-N-1} 
\label{xx3}
\end{equation}
then 
\begin{equation}
| s_a(z)| < |z|^{M+1} | F(z) - a | < |z|^{-N}
\quad \hbox{and} \quad  |s_b(z)| < |z|^{M+1} .
\label{xx4}
\end{equation}

Finally, for any $Q > 0$ there exists $\eta_0 > 0$ such that if
\begin{equation*}
| s_a(z)| + |s_b(z)| \leq \eta_0 
\label{xx5}
\end{equation*}
then
\begin{equation*}
\max \{ | \tan z |, | \cot z | \} \geq Q .
\label{xx6}
\end{equation*}
\end{lem}

\textit{Proof.} 
The definition (\ref{1}) of $F$ gives
\begin{equation}
T - F = \frac{1+T^2}{L+T} , \quad 
s_a = \frac{T(L+T) (F-a) }{1+T^2} =
\frac{T( TL-1 - a(L+T))}{1+T^2} 
\label{C3}
\end{equation}
and hence
$$
s_a(z) = \cos z \sin z ( ( \tan z -a)L(z) -1 - a \tan z ),
$$
which is (\ref{xx2}).

Now suppose that $z$ is large and satisfies (\ref{xx3}). If
\begin{equation}
|T| > c_a = 2|a| + 1
\label{xx7}
\end{equation}
then 
$$
| s_a(z) | = \left| \frac{F(z)-a}{1 - F(z)/T } \right| 
\leq 2| F(z) - a | \quad \hbox{and} \quad
| s_b(z) | = \left| \frac{F(z)-b}{1 - F(z)/T } \right| 
\leq 2| F(z) - b |
$$
and (\ref{xx4}) is obvious, while if (\ref{xx7}) fails then
writing
$$
|T - F(z) | = | T - a - (F(z) - a) | \geq 
\frac12 |z|^{-M} 
$$
leads to
$$
|s_a(z)| \leq 2 c_a |z|^{M} | F(z) - a | 
\quad \hbox{and} \quad |s_b(z)| \leq 2 c_a |z|^{M} | F(z) - b |  ,
$$
which again gives (\ref{xx4}).

Finally, suppose that there exists a sequence
$( z_n ) $ 
such that 
$$
| s_a(z_n)| + | s_b(z_n)| \to 0, \quad
\max \{ | \tan z_n |, | \cot z_n | \} < Q .
$$
Then $|F(z_n)| = O(1) $, because otherwise
writing
$$
T = - \frac{s_a (1-T/F)}{1 - a/F} 
$$
gives a subsequence with
$ \tan z_n = o(1)$, an immediate contradiction. Hence
(\ref{xx1}) yields
$$
|F(z_n) - a| = O( | s_a(z_n) |) = o(1), \quad
|F(z_n) - b| = O( | s_b(z_n) |) =o(1),
$$
which is obviously impossible.
\hfill$\Box$
\vspace{.1in}


\section{Direct transcendental singularities}

Recall the classification of transcendental singularities 
summarised prior to Lemma \ref{directsinglem}.

\begin{lem}\label{lemB1}
If $F^{-1}$ has a direct transcendental singularity 
over $a \in \C \setminus \R$ then $a = \pm i$.
Moreover the function $L$ has finitely many asymptotic values in
$\C \setminus \R$, and $L^{-1}$
cannot have a direct transcendental singularity 
over $a \in \C \setminus \R$.
\end{lem}
\textit{Proof.}
Let $g$ be $F$ or $L$, and  
assume that $g^{-1}$ has a direct transcendental singularity 
over $a \in \C \setminus \R$, with $a \neq \pm i$ if $g = F$. Then there
exist a small
positive $\delta_1$ and a component $D$ of the set
$\{ z \in \C : |g(z) - a| < \delta_1 \}$ such that $g(z) \neq a$ on $D$.
Moreover the function
$$
v(z) = \log \frac{\delta_1}{|g(z) - a|}  \quad (z \in D), \quad
v(z) = 0 \quad (z \in \C \setminus D),
$$
is subharmonic in $\C$. Since $g$ is real meromorphic it may be assumed
that $D \subseteq H$. But $\mathfrak{T}(r, g) = O( \log r )$
as $r \to \infty$ by (\ref{A2}) and (\ref{1}), and so the same argument as
in Lemma \ref{directsinglem} shows that $B(r, v) = O(r \log r)$ as $r
\to \infty$ (compare (\ref{dir0})). In particular $v$ has order
at most $1$.

Let $\delta $ be small and positive and
suppose first that $f \in U_{2p}^*$. If $f$ has at least
one real zero then (\ref{cara}), (\ref{B1}),
(\ref{degphidef}) and (\ref{1}) show that
\begin{equation}
L(z) \to \infty \quad \hbox{and} \quad F(z) \to i  
\quad \hbox{as} \quad z \to \infty, \, \delta < | \arg z | < \pi - \delta .
\label{dir2}
\end{equation}
On the other hand if $f$ has no real zeros then evidently $L$ is
a rational function, with a pole at infinity since $p$ is positive,
and again (\ref{dir2}) holds. Hence
for large $r$ the angular measure of $S(0, r) \cap D$ is at
most $2 \delta $, and by a standard application of Lemma \ref{lemhm}
the order of the subharmonic function $v$ is at least $\pi /2 \delta$.
Since $\delta$ may be chosen arbitrarily small this is a contradiction. 

Suppose next that $f$ has infinite order. Here a different argument is
required since (\ref{dir2}) is not available, and instead a contradiction
will be obtained by
showing that $v$ has lower order greater than $3/2$. The function $\phi$ in 
(\ref{B1}) is transcendental of order at most $1$, by Lemma \ref{lemA2},
and there exists a rational function $R_1$ with at most a simple pole at
infinity such that
$$
\phi_1(z) = \frac{\phi(z) - R_1(z)}{ z^2} 
$$
is entire and transcendental of
order at most $1$. For large $z$ it follows using (\ref{cara})
again that
\begin{equation}
\hbox{if}
\quad | \phi_1(z)| > 1 , \quad 
\delta < | \arg z | < \pi - \delta
\quad \hbox{then} \quad 1/L(z) = o(1),  \quad 
| g(z) - a | \geq \delta_1.
\label{dir3}
\end{equation}
Let $C$ be a component of the set $\{ z \in \C : |\phi_1(z)| > 1 \}$
and for $s > 0$ let $\theta_C (s), \theta_D (s)$ denote the
angular measure of $C \cap S(0, s), D \cap S(0, s)$ respectively.
Since $g^{-1}$ also has a direct transcendental singularity 
over $\bar a $, it follows from (\ref{dir3}) that, for large $s$,
\begin{equation}
\theta_C(s) + 2 \theta_D(s) \leq 2 \pi + 4 \delta .
\label{dir4}
\end{equation}
Let $\theta_C^*(s) = \infty$ if $S(0, s) \subseteq C$ and
$\theta_C^*(s) = \theta_C(s)$ otherwise. Then
\begin{equation}
9 \leq \left(
\frac1{\theta_C^*(s)} + \frac2{\theta_D(s)} \right)(2 \pi + 4 \delta )
\label{dir5}
\end{equation}
for large $s$, using (\ref{dir3}), (\ref{dir4}), 
the Cauchy-Schwarz inequality
and the fact that $\delta$ is small.
Integrating (\ref{dir5}) from $r_0 $ to $r$, where $r_0$ is large,
and using Lemma \ref{lemhm} yields, as $r \to \infty$,
\begin{equation}
\frac{9 \log r}{2 + 4 \delta /\pi } - O(1) \leq
\int_{r_0}^r 
\left(
\frac1{\theta_C^*(s)} + \frac2{\theta_D(s)} \right) \, \frac{\pi \, ds}{s}
\leq (1 + o(1)) \log r +
2 \int_{r_0}^r \frac{\pi \, ds}{s \theta_D(s)},
\label{dir6}
\end{equation}
since $\phi_1$ has order at
most $1$. Applying
Lemma \ref{lemhm} to $v$, and using 
(\ref{dir6}) and the fact that $\delta $ is small by assumption
now shows that the lower order of $v$ is
at least
$$
\frac12 \left( \frac{9 }{2 + 4 \delta /\pi} - 1 \right) > \frac32 .
$$

Finally, the assertion
that $L$ cannot have infinitely many asymptotic values 
$a \in \C \setminus \R$ is proved by observing that in the contrary case
$L^{-1}$ would have at least two direct transcendental singularities
over $\infty$ lying in $H$, which by (\ref{A2}) contradicts 
Lemma \ref{directsinglem}.
\hfill$\Box$
\vspace{.1in}

\section{Indirect transcendental singularities}

The following proposition uses again the terminology summarised prior
to Lemma \ref{directsinglem}.

\begin{prop}\label{prop1}
There does not exist $\alpha \in \C \setminus \R $ such that
the inverse function $F^{-1}$ has
an indirect transcendental singularity over $\alpha$. 
\end{prop} 
\textit{Proof.} 
To establish this proposition will require the whole of this
section and a number of intermediate lemmas.
Assume that there exists $\alpha \in \C \setminus \R$
such that $F^{-1}$ has
an indirect transcendental singularity over $\alpha$. Since $F$ is real it may
be assumed that the corresponding path and components lie in $H$. 
The key idea 
will be to show that
there exist paths $\Gamma_j$ tending to infinity in $H$ on which 
$F(z)$ tends to distinct values $\beta_j \in H$, and 
to use the fact that for most large $z$
on $\Gamma_j$ it follows from Lemma \ref{lemC3} that the function 
$s_{\beta_j}(z)$, as defined
by (\ref{xx1}), is small. A contradiction will then arise from an argument of
Phragm\'en-Lindel\"of type, using the fact that 
these functions $s_{\beta_j}(z)$ have finitely many
poles in $H$. Unfortunately, however, complications arise because in
principle the $\Gamma_j$ may pass close to points where
$\tan z = \beta_j$, and near these points Lemma \ref{lemC3} cannot be applied.

\begin{lem}\label{lem3}
Let $N  $ be a large positive integer.
There exist pairwise distinct complex numbers 
$$\beta_j \in \C \setminus \R , \quad
\beta_j \neq \pm i  , \quad j = 1, 2, \ldots, N , $$
with $| \beta_j - \alpha | = \eta_j $ small and positive,
and pairwise disjoint simply connected
domains $U_j \subseteq H$ with the following properties:\\
(i) F maps $U_j$ univalently onto the disc $D(\alpha, \eta_j)$;\\
(ii) there exists a simple path $\Gamma_j$ tending to infinity
in $U_j$,
mapped by $F$ onto the half-open line segment $[\alpha, \beta_j)$,
such that $F(z) \to \beta_j$ as $z$ tends to infinity on $\Gamma_j$.\\
\end{lem}
\textit{Proof.} The existence of $\beta_j$, $U_j$ and $\Gamma_j$
with $0 \neq | \beta_j - \alpha | \to 0$ as $j \to \infty$ 
follows from the definition \cite{BE} of an indirect singularity
(see also \cite[Lemma 10.3, p.370]{lajda}), and in particular $\beta_j \neq 
\pm i$
for $j$ sufficiently large.
\hfill$\Box$
\vspace{.1in}

\begin{lem}\label{lemF0}
Let $r$ be large and positive. Then the domains $U_j$ and paths $\Gamma_j$ of
Lemma \ref{lem3} may be labelled so that
\begin{equation}
\int_{r^{1/32}}^{r^{1/16}} \frac{ \pi  dt}{t \theta_{U_j} (t)} 
> 2048 \pi \log r 
\label{F1}
\end{equation}
for $j = 1, \ldots, 500$, where $\theta_{U_j} (t)$ is as defined in
Lemma \ref{lemE1}. 
\end{lem}
\textit{Proof.} 
Lemma \ref{lemE1} shows that (\ref{F1}) can fail
for at most $2 \cdot 2048 \pi \cdot 32  $ of the domains $U_j$, and the result
follows since $N$ is large.
\hfill$\Box$
\vspace{.1in}

For the remainder of this section $r$ as in Lemma \ref{lemF0} will be fixed,
and $d , d'$ will be used to denote positive constants, not necessarily the
same at each occurrence, but always independent of $r$.
\begin{lem}\label{lemF1}
The function $F$ satisfies, for $j = 1, \ldots, 500$,
\begin{equation}
|F(z) - \beta_j | \leq d r^{-512} \leq d |z|^{-64}  \quad
\hbox{for $z \in \Gamma_j, \, r^{1/8} \leq |z| \leq r^{8} $.}
\label{F2}
\end{equation}
\end{lem}
\textit{Proof.} 
This uses the argument of \cite[p.371]{lajda}. 
Let $G$ be that branch of the inverse function $F^{-1}$ mapping
$D(\alpha, \eta_j )$ onto $U_j$. For 
$u \in \Gamma_j$ the distance from $u$ to
$\partial U$ is at most 
$|u| \theta_{U_j} (|u|)$ and so Koebe's theorem 
implies that
$$
|(v - \beta_j ) G'(v)| \leq 4 |u| \theta_{U_j}(|u|) 
\quad \hbox{for} \quad u = G(v), \, v \in [\alpha , \beta_j) .
$$
Hence, for $z \in \Gamma_j$ with $r^{1/8} \leq |z| \leq r^{8}$
writing $w = F(z)$ and
$u = G(v)$ for $v \in [\alpha, w]$ gives, using (\ref{F1}), 
\begin{eqnarray*}
\log \left| \frac{\beta_j-\alpha}{\beta_j - F(z)}  \right|
&=& \int_\alpha^w \frac{|dv|}{|\beta_j-v|}
= \int_{G(\alpha)}^z \frac{|du|}{|(\beta_j-v)G'(v)|}
\geq \int_{G(\alpha)}^z \frac{|du|}{4 |u| \theta_{U_j}(|u|)}\\
&\geq& \int_{r^{1/32}}^{r^{1/16}} \frac{  dt}{4 t \theta_{U_j} (t)} 
> 512 \log r \geq 64 \log |z|. 
\end{eqnarray*} 
\hfill$\Box$
\vspace{.1in}

\begin{lem}\label{lemF2}
For $1 \leq j \leq 500$ pick an arc $\lambda_j$ of $\Gamma_j$ joining
$S(0, r^{1/8} )$ to $S(0, r^{8} )$ and, apart from its endpoints,
lying in $A(r^{1/8}, r^{8})$. 
By re-labelling if necessary it may be assumed
that these arcs $\lambda_j$ separate
the half-annulus $A^+(r^{1/8}, r^{8} )$ in counter-clockwise
order.
For $1 \leq j < 500$ let $W_j$ be the part of $A^+(r^{1/8}, r^{8})$
separating $\lambda_j$ from $\lambda_{j+1}$. Then there exists 
$k $ such that
\begin{equation}
\int_{r^{1/4}}^{r^{1/3}} \frac{\pi \, dt}{t\theta_{W_k} (t)}
> 16 \log r 
\quad \hbox{and} \quad
\int_{r^3}^{r^{4}} \frac{\pi \, dt}{t\theta_{W_k} (t)} > 16 \log r .
\label{F15}
\end{equation}
\end{lem}

\textit{Proof.} The existence of $k$ satisfying (\ref{F15}) follows since
Lemma \ref{lemE1} shows that the first inequality of (\ref{F15})
fails for at most $2 \cdot 16 \cdot 12
= 384$ of the $W_j$ and the second for at most
$ 32  $ of them. 
\hfill$\Box$
\vspace{.1in}

\begin{lem}\label{lemF5}
Choose $k$ satisfying $(\ref{F15})$ and for convenience write
\begin{equation*}
a = \beta_k, \quad b = \beta_{k+1}, \quad
\lambda_a = \lambda_k , \quad \lambda_b = \lambda_{k+1}. 
\end{equation*}
Denote by $u_\nu$ the solutions in the annulus
$A(r^{1/16}, r^{16} )$ of the equations
\begin{equation}
\tan z = a, b     .
\label{F11}
\end{equation}
Then the discs $D(u_\nu, |u_\nu|^{-2} )$ are pairwise disjoint. 
Next, set
\begin{equation}
P(z) = s_a(z) s_b(z) .
\label{F17}
\end{equation}
Then $P$ has no poles in $A^+(r^{1/16}, r^{16})$. Finally, the functions
$s_a, s_b$ and $P$ 
satisfy the estimates:\\
(i)
\begin{equation}
|s_a(z)   | \leq |z|^{-28} \quad 
\hbox{and} \quad |P(z)| \leq |z|^{-14} \quad
\hbox{for all $z \in \lambda_a  \setminus \bigcup_\nu D(u_\nu, 
|u_\nu|^{-12} )$};
\label{F18}
\end{equation}
(ii) 
\begin{equation}
|s_b(z)   | \leq |z|^{-28} \quad 
\hbox{and} \quad |P(z)| \leq |z|^{-14} \quad
\hbox{for all $z \in \lambda_b  \setminus \bigcup_\nu D(u_\nu, 
|u_\nu|^{-12} )$};
\label{F19}
\end{equation}
(iii) 
\begin{equation}
|P(z)| \leq \exp( d |u_\nu| \log |u_\nu| ) 
 \quad
\hbox{for all $z \in  D(u_\nu, |u_\nu|^{-12} )$}.
\label{F20}
\end{equation}
\end{lem}
\textit{Proof.} 
The assertion concerning the discs $D(u_\nu, |u_\nu|^{-2} )$ is
obvious since $a \neq b$ and $r$ is large, and $P$ has no poles
in $A^+(r^{1/16}, r^{16})$ by (\ref{xx2}).
To prove (i) let 
$z \in \lambda_a  \setminus \bigcup_\nu D(u_\nu, 
|u_\nu|^{-12} )$ and observe first that
the choice of the $u_\nu$ gives
$ | \tan z - a | \geq |z|^{-13} $. 
Recalling (\ref{F2}) and applying
Lemma \ref{lemC3} with $M = 13, N = 28$ now leads at once to
(\ref{F18}), and
(\ref{F19}) is obtained using the same argument.

Finally, to prove (iii)
suppose that $z \in  D(u_\nu, |u_\nu|^{-12} )$. Then
$d \leq | {\rm Im} \, z | \leq d'$, 
since $a$ and $b$ are non-real. Hence (\ref{cara}) 
and (\ref{B1}) give
$$
| \psi (z)| \leq d |z|^2, \quad 
|L(z)| \leq d |z|^2 M( |z|, \phi ),
$$
and 
(\ref{F20}) follows using (\ref{phifest})
and (\ref{xx2}). This proves Lemma \ref{lemF5}.
\hfill$\Box$
\vspace{.1in}

\begin{lem}\label{lemF6}
Let $k, a, b, \lambda_a, \lambda_b$
be as in Lemma \ref{lemF5}. Then there exist
$R_0, R_1, R_2$ satisfying 
\begin{equation}
R_0 , R_1, R_2 \not \in E_0, \quad
r^{1/2} \leq R_0 \leq r^2, \quad
r^{1/7} \leq R_1 \leq  r^{1/6} , \quad
r^{6} \leq R_2 \leq r^{7},   
\label{G1}
\end{equation}
where $E_0$ is the exceptional set of Lemma \ref{lemA2}, and 
with the additional properties that
\begin{equation}
S(0, R_\mu) \cap D(u_\nu, |u_\nu|^{-2}) = \emptyset
\label{G2}
\end{equation}
for $\mu = 0, 1, 2$ and each $\nu $, as well as
\begin{equation}
| \tan z | + | \cot z | \leq d \quad
\hbox{for} \quad |z| = R_0 .
\label{xx11}
\end{equation}
Finally, there exists
$w_k \in W_k \cap S(0, R_0)$ with
\begin{equation}
|s_a(w_{k})| = |s_b(w_{k})|.
\label{G3}
\end{equation}
\end{lem}
\textit{Proof.}
First, $R_0, R_1$ and $R_2$ exist because $E_0$ has finite logarithmic measure
and $r$ is large, and  
the discs $D(u_\nu, |u_\nu|^{-2})$ have sum of radii at most $d$. 
To prove the existence of $w_k$ observe that since $W_k$ 
separates $\lambda_a$ from $\lambda_b$ there exists an
arc $A_k $ of the circle $S(0, R_0)$ which lies in $W_k$ apart from its
endpoints $v_a$ and $v_{b}$, which satisfy 
$v_a \in \lambda_a $ and $v_b \in \lambda_b$. 
It then follows using Lemma \ref{lemC3},
(\ref{F18}) and (\ref{F19})  
that
$$
|s_b(v_a) | > | s_a(v_a)| , \quad
|s_a(v_b) | > | s_b(v_b)|,
$$
and so a point
$w_k \in A_k$ satisfying (\ref{G3}) exists by continuity. 
\hfill$\Box$
\vspace{.1in}

A contradiction will now be obtained using
harmonic measure. 
Let $k$, $W_k$ and $w_k$ be as in Lemmas \ref{lemF5} and \ref{lemF6},
and let $D$ be the component of the set
$W_k \cap A^+(R_1, R_2)$ which contains $w_k$.
The function 
\begin{equation}
Q(z) = z^{14} P(z) ,
\label{G4}
\end{equation}
is analytic on the closure of $D$,
and evidently
\begin{equation}
|Q(z)| \leq 1  \quad \hbox{for all $z \in (\lambda_a \cup \lambda_{b})
\setminus \bigcup_\nu D( u_\nu , |u_\nu|^{-12} )$}
\label{G5}
\end{equation}
by (\ref{F18}) and (\ref{F19}).
Next,
\begin{equation}
|Q(z)| \leq 
R_\mu^{14} \exp( d R_\mu \log R_\mu) \leq \exp( d r^{8} \log r ) \quad 
\hbox{for all $z \in S(0, R_\mu)$, $\mu = 1,2$,}
\label{G6}
\end{equation}
by Lemma \ref{lemA2}, (\ref{xx2}) and (\ref{G1}), while (\ref{F15}) 
and (\ref{G1}) give a harmonic measure estimate
\begin{equation}
\omega(w_k, D, S(0, R_\mu)  \cap \partial D ) \leq 
d r^{-16}  \quad \hbox{for} \quad \mu = 1, 2.
\label{G7}
\end{equation}

It remains to consider the intersection of 
$\partial D$ with the discs $D( u_\nu , |u_\nu|^{-12} )$.
First, (\ref{F20}) and (\ref{G4}) yield
\begin{equation}
\log |Q(z)| \leq 
d |u_\nu| \log |u_\nu| \quad
\hbox{for all $z \in  D(u_\nu, |u_\nu|^{-12} )$}.
\label{meet2}
\end{equation}
Suppose then that $D( u_\nu , |u_\nu|^{-12} )$ meets $\partial D$,
at $y_\nu$ say. Then it follows that
\begin{equation}
S(y_\nu, t ) \setminus D \neq \emptyset \quad \hbox{for} \quad
|u_\nu|^{-11} \leq t \leq |u_\nu|^{-3} .
\label{meet}
\end{equation}
For if such a circle
$S(y_\nu, t )$ lies in $D$ then the closed disc $E^*$ given by
$|z - y_\nu | \leq t$ lies in each of the simply connected domains
$A^+(R_1, R_2)$ and $W_k$, and 
so in some component of the intersection; but then $E^* \subseteq D$, since
$E^*$ meets $D$ near $y_\nu$, which contradicts the fact that 
$y_\nu \not \in D$ and proves (\ref{meet}).
But   
$$
|w_k - y_\nu | \geq |w_k - u_\nu | -
| u_\nu - y_\nu  |  \geq |u_\nu |^{-2} - |u_\nu|^{-12} \geq  |u_\nu|^{-3},
$$
since 
$w_k \in S(0, R_0)$. 
It now follows that
the change of variables
$$\zeta = \frac1{z - y_\nu}, \quad
\zeta_k = \zeta( w_k) = \frac1{w_k  - y_\nu}, 
$$
maps $D$ to a domain $D^*$ in $\C$ such that 
the circle $S(0, t)$ meets $\C \setminus D^*$ for 
$|u_\nu|^3 \leq t \leq |u_\nu|^{11}$, while
$$| \zeta_k| \leq |u_\nu|^3 \quad \hbox{and} \quad 
| \zeta (z) | \geq  |u_\nu|^{11} \quad \hbox{for} \quad z \in 
D( u_\nu , |u_\nu|^{-12} ).$$
This now gives, by conformal invariance of harmonic measure,
\begin{eqnarray}
\omega(w_k, D,  D( u_\nu , |u_\nu|^{-12} ) \cap \partial D)
&\leq& d \exp \left( - \int_{|u_k|^4}^{|u_k|^{10}} \frac{\pi \, dt}{t
\theta_{D^*}(t)} \right) \nonumber \\
&\leq& d \exp \left( - \int_{|u_k|^4}^{|u_k|^{10}} \frac{ dt}{2t} \right)
 \leq d |u_\nu|^{-3} .
\label{G8}
\end{eqnarray}

Combining (\ref{G5}), (\ref{G6}), (\ref{G7}), 
(\ref{meet2}) and (\ref{G8}) leads to
$$
\log |Q(w_k)| \leq d \left(  r^{8 - 16 } \log r +
\sum_\nu |u_\nu |^{-2} \log |u_\nu| \right) \leq d.
$$
Using (\ref{F17}), (\ref{G3}) and (\ref{G4}) it now follows that
$$s_a(w_k) = o(1), \quad s_b(w_k) = o(1),$$
which contradicts Lemma \ref{lemC3}
and proves Proposition \ref{prop1}.

\section{Zeros of $\phi$}

\begin{lem}\label{phizeroslem}
If $f$ has infinite order then the function $\phi$ has infinitely many
zeros.
\end{lem}
\textit{Proof.}
Assume that $f$ has infinite order but $\phi$ has finitely many
zeros. Then it follows from Lemma \ref{factorlem}(iii) and Lemma \ref{lemA2}
that there exist a rational function $R_1$ and a non-zero real constant 
$c_1$ such that
\begin{equation}
\phi(z) = R_1(z) e^{c_1 z} .
\label{dir10}
\end{equation}
Hence it follows using (\ref{cara}), (\ref{B1}) and (\ref{1}) that 
$F(z) \to i$ as $z \to \infty$ on each of the rays 
$L_1, L_2$ given by $\arg z = \pi /2 \pm \pi /16$. Thus each of the rays
$L_1, L_2$ gives rise to a transcendental singularity of 
$F^{-1}$ over $i$, which must be direct by Proposition \ref{prop1}. 
Applying (\ref{A2}) and (\ref{1}) in combination with Lemma
\ref{directsinglem} then shows that the two rays $L_1, L_2$ must determine
the same direct transcendental singularity of 
$F^{-1}$, and so there exist a small positive constant 
$\delta$ and a component
$C$ of the set $\{ z \in \C : |F(z) - i | < \delta \}$, on which
$F(z) \neq i$, such that $z \in C$ for all large $z$ with
$\arg z = \pi /2 \pm \pi /16$. It follows from (\ref{L=i}) that
$L(z) \neq i$ on $C$.

Since $i$ is not a limit point of transcendental singularities of 
$F^{-1}$, by Lemma
\ref{directsinglem} and Proposition \ref{prop1}, nor of critical values of
$F$, by Lemma \ref{lem1}, the singularity over $i$ is logarithmic.
Provided $\delta $ is small enough this implies  
in particular that the boundary of $C$ consists of one simple curve tending
to infinity in both directions \cite{Nev}
(see also \S\ref{rh}). Hence all large $z$ with 
$| \arg z - \pi /2 | \leq \pi /16 $ are in $C$. But it is evident
from (\ref{cara}), (\ref{B1}),
(\ref{dir10}), the Phragm\'en-Lindel\"of principle
and the fact that $c_1$ is real
that the equation $L(z) = i$ must have infinitely
many solutions near the positive imaginary axis, and this is a contradiction. 
\hfill$\Box$
\vspace{.1in}

\section{The behaviour of $L$ near zeros of $\phi$}
The next lemma uses the notation of Definition \ref{def1}, the 
stated convention that
all counts of zeros are with regard to multiplicity unless indicated
otherwise, and
reasoning similar to \cite[p.984]{BEL} and \cite[Lemma 14.1]{lajda}.

\begin{lem}\label{componentslem}
If $f \in U_{2p}^*$ and $f$ has $2q$ distinct non-real zeros then 
for sufficiently small positive $\lambda $ there
are at least $p + q$ bounded components $K_j \subseteq H$ of the set
$L^{-1} ( D^+(0, \lambda ))$, each mapped univalently onto 
$D^+(0, \lambda )$ by $L$, and with a zero of $L$ on $\partial K_j$.
If $f$ has infinite order and $M \in \N$ then for
sufficiently small positive $\lambda $ there exist at least $M$
such components $K_j$.
\end{lem}
\textit{Proof.}
Suppose first that
$\zeta \in H$ is a zero of $L$ of multiplicity $m$. Then since 
$ L_\zeta (z) = L(z)^{1/m} $ is analytic and
univalent near $\zeta$, it follows that provided $\lambda$ is small
enough there exist $m$ components $K_j$ as in the statement of the lemma.

Next, let $\zeta$ be a real zero of $L$ of even multiplicity $m$.
Then for $\lambda$ small enough $\zeta$ gives rise
to $m/2$ components $K_j$, and as $x$ passes through $\zeta$ from
left to right the sign of $L(x)$ does not change.
Now suppose that $\zeta$ is a real zero of $L$ of odd multiplicity $m$
and that $\lambda$ is small. If $L^{(m)}(\zeta) > 0$ then 
$\zeta$ gives rise
to $(m+1)/2$ components $K_j$, and $L(x)$
has a positive sign change at $\zeta$, that is,
as $x$ passes through $\zeta$ from
left to right the sign of $L(x)$ changes from negative to positive.
On the other hand if $L^{(m)}(\zeta) < 0$ then $\zeta$ gives rise
to $(m-1)/2$ components $K_j$, and
$L(x)$ has a negative sign change at $\zeta$.

In the case where $f$ has infinite order, it is 
now clear that the conclusion of the lemma holds if $L$ has infinitely
many non-real or multiple zeros, so assume that all but finitely many zeros
of $L$ are real and simple. Since $\phi$ has infinitely
many zeros by Lemma \ref{phizeroslem}, there are two
alternatives. The first is that there exists
an unbounded open interval $I$ of $\R$ containing no poles of $L$ but
infinitely many zeros of $\phi$ and so of $L$, in which case
$I$ evidently contains infinitely many zeros $\zeta$ of $L$ with
$L'(\zeta ) > 0$, and the conclusion of the lemma follows.
The second alternative is that there exist
infinitely many bounded open intervals
$I = (a, b)$ lying between adjacent zeros $a, b$ of $f$ and containing
at least one zero of $\phi$, in which case $\psi$ has a zero in $(a, b)$ by
construction, or by the fact that $\psi$ has negative residues,
and so $L$ has at least two zeros $\zeta \in (a, b)$, at
least one of them having $L'(\zeta ) \geq 0$, so that again the conclusion of 
the lemma follows.

Suppose now that $f \in U_{2p}^*$. Then $\phi$ is a rational function.
Let $I$ be a component of $\R \setminus f^{-1}( \{ 0 \} )$ containing
$\mu_I > 0$ zeros of $\phi$ and
$m_I $ zeros of $L$. Then $m_I \geq \mu_I$ and $\mu_I$ is even
by Lemma \ref{factorlem}(iv).
Hence, by the above analysis, if $\lambda $ is sufficiently small,
the interval $I$ gives rise to
\begin{equation}
n_I = \frac{m_I + s_I}2 \geq \frac{\mu_I + s_I}2
\label{zeros1}
\end{equation}
components $K_j$, where 
$s_I$ is the number of
positive sign changes minus the number of negative sign changes 
undergone by $L(x)$ on $I$. Since $s_I \geq -1$ and $\mu_I$ is
even, (\ref{zeros1}) yields $n_I \geq \mu_I /2$. 

Let $2r$ be the number of non-real zeros of $\phi$, these coinciding with
zeros of $L$, and let $2 \nu $ be the number of real zeros of $\phi$.
Summing over all the intervals $I$ 
it follows that for small enough $\lambda$
there are at least
$\nu + r $ components $K_j$ as in the statement of the lemma.
But each non-real zero of $f$
is a simple pole of $\phi$, and the argument principle gives
\begin{equation}
2 \nu + 2r = 2q + {\rm deg}_\infty (\phi) .
\label{zeros2}
\end{equation}
Thus the conclusion of the lemma follows at once from (\ref{degphidef}),
except in the case where $f$ has no real zeros. In this last case, however,
$L$ is a rational function and $f$ satisfies
$$
f = P_0 \exp( P_1) 
$$
where $P_0$ is a real polynomial with no real zeros, and $P_1$ is
a real polynomial of degree $d_1 \geq 2p$. If $d_1 \geq 2p+1$ then
${\rm deg}_\infty (\phi) \geq 2p-1$ by
Lemma \ref{factorlem}(i) and (\ref{cara}),
and again the result follows from (\ref{zeros2}). Suppose finally
that $d_1 = 2p$. Then the leading coefficient $c_1$ of $P_1$ is positive
and $L(z) \sim 2p c_1 z^{2p-1}$ as $z \to \infty$. Here
there is one component $I = \R$, and $s_I = 1$, and if $\lambda$ is small
then
(\ref{zeros1}) and the argument principle applied to $L$
give at least $p^*$ components $K_j$, where
$$
p^* \geq n_I + r = \frac{m_I + 1}2 + r =
\frac{m_I + 2r + 1}2 = \frac{2p-1 + 2q +1}2 = p + q .
$$
\hfill$\Box$
\vspace{.1in}

\section{Components of $W$ }

This section will discuss components of
the set $W$ defined in (\ref{2}). Recall from Lemma \ref{lem1} that $F$
has finitely many non-real critical points, and that by Lemma \ref{lemB1}
and Proposition \ref{prop1} the only possible asymptotic value $w \in H$
of $F$ is $i$.

\begin{lem}\label{Wlem}
Choose a simple polygonal path $\Lambda $ in 
$(H \cup \{ 0 \}) \setminus \{ i \}$, such that $\Lambda$ contains all
critical values of $F$ in $H \setminus \{ i \}$, and let $H^* =
H \setminus \Lambda $. Then all components $A^*$ of the set
$W^* = \{ z \in H: F(z) \in H^* \}$ are
simply connected. Moreover each such component $A^*$ belongs to one of two
types:\\
(a) type I, for which $A^*$ contains no $i$-points of $F$, but a path  
tending to infinity on which $F(z) \to i$, and $A^*$ is mapped onto
$H^* \setminus \{ i \}$ by $F$;\\
(b) type II, for which $A^*$ contains one $i$-point of $F$, of multiplicity
$m$, and is mapped $m: 1$ onto $H^*$ by $F$.

There is at most one type I component $A^*$ of $W^*$, and the following
properties hold:\\
(i) each
component $A$ of $W$ contains finitely many components $A^*$ of 
$W^*$ and so finitely many $i$-points of $F$;\\
(ii) if a component $A$ of $W$ does not contain a
type I component $A^*$ of $W^*$ and does not
contain any critical points of $F$ then $A$ is mapped conformally onto
$H$ by $F$. 
\end{lem}
\textit{Proof.}
Applying the standard transformation
\begin{equation}
u = G(w) = \frac{F(w)-i}{F(w)+i} 
\label{uGdef}
\end{equation}
and recalling that $F(\R) \subseteq \R \cup \{ \infty \}$ shows that
every component of $W$ is a component of the set
$\{ w \in \C :  |G(w)|  < 1 \}$. 
Hence the fact that the components $A^*$ are simply
connected, their classification as types I or II, and
properties (i) and (ii) all follow from the discussion in \S\ref{rh}.
Since every type I component of $W^*$ gives rise to a direct singularity
of $F^{-1}$ lying in $H$, it follows from Lemma \ref{directsinglem},
(\ref{A2}) and (\ref{1}) that
there is at most one type I component. 
\hfill$\Box$
\vspace{.1in}

\begin{lem}\label{Wlem2}
Let $A$ be a component of $W$ containing a type I component of
$W^*$. Then the number of $i$-points of $F$ in $A$ is at most the
number of zeros of $F'$ in $A$.
\end{lem}
\textit{Proof.}
This follows from Lemmas \ref{rhlem} and \ref{Wlem}.
\hfill$\Box$
\vspace{.1in}

\section{Components of $Y$}\label{Y}

Recall from Lemma \ref{componentslem} that there exist a small
positive $\lambda $ and $M$ components
$K_j \subseteq H$ of the set $L^{-1}( D^+(0, \lambda ) )$, each mapped
univalently onto $D^+(0, \lambda )$ by $L$. Here $M = p+q$ if
$f \in U_{2p}^*$, where $2q$ is the number of distinct non-real zeros
of $f$, and $M$ may be chosen arbitrarily large if $f$ has
infinite order. 

Each such $K_j$ lies in a component $C_j$ of the set
$Y$ defined in (\ref{2}), which in turn lies
in a component $A_j$ of $W$, by Lemma \ref{lem1}. Here the $C_j$ 
corresponding to different $K_j$ need not be distinct, and this is
also the case for the $A_j$ corresponding to different $C_j$.

\begin{lem}\label{Cjlem2}
For each $C_\nu$ the number of $K_j$ contained in $C_\nu$ is
at most the number of $i$-points of $F$ in $C_\nu$, and this number
is finite.
\end{lem}
\textit{Proof.}
Recall first that the number of $i$-points of $L$ in $C_\nu$ 
equals the number of $i$-points of $F$ in $C_\nu$, by
(\ref{L=i}), and since $C_\nu$ lies in some $A_\mu$ this number is finite, by
Lemma \ref{Wlem}. Choose a circular arc $\gamma$ joining $0$ to $i$ in
the closure of $H$ and passing through no singular values of $L^{-1}$
apart possibly from $0$ and $i$ themselves. This is possible
by Lemma \ref{lemB1}. For each $K_j \subseteq C_\nu$ 
choose $z_j \in K_j$ with $L(z_j) \in \gamma$. Then 
the inverse function 
$L^{-1}$ may be continued along the half-open 
subarc of $\gamma$  joining $L(z_j)$ to
$i$, by the choice of $\gamma$, and the image $\gamma_j (w)$ 
of this continuation starts at $z_j$ and lies in $C_\nu$.
If $\gamma_j (w)$ tends
to infinity as $w \to i$ this gives a path tending to infinity in $C_\nu$
on which $L(z)$
tends to $i$. But an indirect singularity
of $L^{-1}$ over $i$ is excluded since there are finitely many
$i$-points of $L$ in $C_\nu$, while a direct singularity is ruled out
by Lemma \ref{lemB1}. 

Hence $\gamma_j (w)$ cannot tend to infinity,
so that $\gamma_j(w)$  has a finite 
limit point $z_j^*$ as $w$ tends to $i$ along $\gamma$. Thus 
$\gamma_j (w)$ tends to $z_j^*$, and
$z_j^*$ must be an $i$-point of $L$ in $C_\nu$. Moreover the
number of such $\gamma_j$ tending to an $i$-point of $L$ in $C_\nu$
is at most the multiplicity of that $i$-point, 
which is the same for $L$ as for $F$,
by (\ref{L=i}). This proves the lemma.
\hfill$\Box$
\vspace{.1in}

Now choose $\theta' \in (\pi /4, 3 \pi /4)$ such that
the ray $\gamma'$ given by $z = s e^{i \theta'}$, $0 < s < \infty$,
contains no singular values of
$L^{-1}$, again using Lemma \ref{lemB1}. 
For each $K_j$ choose $z_j' \in K_j$ with $L(z_j') \in \gamma'$,
and continue $L^{-1}$ along $\gamma'$ in the direction of 
$\infty$. Let $\Gamma_j$ be the image of this continuation starting
at $z_j'$. Then $\Gamma_j$ is a path in $C_j$ on which $L(z) \to \infty$,
and $\Gamma_j$ tends either to infinity or to a pole of $L$, which
must be a zero
of $f$ in $H$, by Lemma \ref{lem1}. 
A component $A_\nu$ of $W$ will be called type $(\alpha)$ if
there exists $K_j \subseteq C_j \subseteq A_\nu$ such that
$\Gamma_j$ tends to infinity, and type $(\beta)$ otherwise.

\begin{lem}\label{Cjlem3}
Let $A_\nu$ be type $(\beta)$. Then the number of $K_j$ contained
in $A_\nu$ is at most the number of distinct non-real zeros of
$f$ in $A_\nu$.
\end{lem} 
{\textit{Proof.} For each $K_j$ contained in $A_\nu$ the path $\Gamma_j$
must tend to a zero $v_j$ of $f$ in $H$, and since these are simple
poles of $L$ the $v_j$ for different $K_j$ must be distinct. Moreover,
(\ref{1}) gives $F(v_j) = \tan v_j \in H$ and so $v_j \in A_\nu$.
\hfill$\Box$
\vspace{.1in}

\section{Completion of the proof when $f$ has finite order}
\label{finiteorder}

\begin{lem}\label{Cjlem4}
Assume that $f \in U_{2p}^*$ and let $A_\nu$ be a type $(\alpha)$
component of $W$. 
Then the number of $K_j$ contained
in $A_\nu$ is at most the number of distinct non-real zeros of
$f$ in $A_\nu$ plus the number of zeros of $f''+f$ in $A_\nu$.
\end{lem} 
{\textit{Proof.} 
By Lemmas \ref{lem1} and \ref{Cjlem2} it suffices to show that the number
of $i$-points of $F$ in $A_\nu$ is at most the number of zeros of
$F'$ in $A_\nu$. This follows in turn from Lemma \ref{Wlem2} provided that
it can be shown that $A_\nu$ contains a type I component of the  
set $W^*$ defined in Lemma \ref{Wlem}.

Let the type II components of $W^*$ which are contained in $A_\nu$ be
$A_1^*, \ldots , A_\mu^*$, and let $\eta$ be small and positive. 
Then since each $A_j^*$ is mapped $m_j:1$ onto $H^*$ by $F$, 
for some integer $m_j$,
each set
$$
B_j^* = \{ z \in A_j^* : | F(z) - i | < \eta \} 
$$
is bounded. It suffices therefore to show that there exist points
$z \in A_\nu$ with $|z| $ arbitrarily large and $| F(z) - i | < \eta$,
since these points $z$ must then lie in a type I component of $W^*$.

Let $C$ be a component of $Y$ with $C \subseteq A_\nu$ such that 
$C$ contains a curve $\Gamma_j$ as defined in
\S\ref{Y} which tends to infinity.
Such a component $C$ exists since $A_\nu$ is type $(\alpha)$. 
Choose $R^* , S^* \in (0, \infty)$ such that all non-real zeros of $f$ lie in
$D(0, R^*)$ and $|L(z)| \leq S^*$ on $S(0, R^*)$. Then $C$ contains
an unbounded component $C^*$ of the set
$\{ z \in \C : {\rm Im} \, L(z) > 2 S^* \}$ with no poles of $L$ 
in its closure, using Lemma \ref{lem1}.
The function $v_C$ defined in analogy with (\ref{4})
by
$$
v_C(z) = {\rm Im} \, L(z) \quad (z \in C^*), \quad
v_C(z) = 2 S^* \quad (z \not \in C^*),
$$
is non-constant and subharmonic
in the plane, and of lower order at least $1$ since $v_C = 2S^*$ on 
$\C \setminus H$.
On the other hand $v_C$ has finite order since $L = f'/f$ and $f$ has
finite order. Hence combining Lemma \ref{lemhm}
with a result of Hayman
\cite{Hay3} shows that there exist
positive constants $d_1, d_2, d_3$ and arbitrarily large positive $r$
such that 
$$
B(r, v_C) \leq 3 T(2r, v_C) \leq d_1 T(r, v_C), \quad
v_C(z) > d_2 T(r, v_C) > r^{1-o(1)} 
$$
on a subset of $S(0, r)$ of angular measure at least $d_3$. Therefore
choosing such
points $z$ with $d_3/4 \leq \arg z \leq \pi - d_3 /4$ gives 
$F(z) \sim i$ by (\ref{1}), as required.
\hfill$\Box$
\vspace{.1in}

This completes the proof of Theorem \ref{thm1} when $f \in U_{2p}^*$,
since Lemma \ref{componentslem} gives $p+q$ components $K_j$, but by 
Lemmas \ref{Cjlem3} and \ref{Cjlem4} the number of $K_j$ does not exceed
the number $q$ of distinct zeros of $f$ in $H$ plus the
number of zeros of $f''+f$ in $H$.

\section{Completion of the proof when $f$ has infinite order}

Assume now that $f$ has infinite order. Here the method of Lemma \ref{Cjlem4}
is not available, and a different approach
is required, based on the notation and results of \S\ref{Y}. 

\begin{lem}\label{lemG3}
Let $N$ be a positive integer. 
Then there exist at least $N$ distinct
components $A$ of $W$ with the following
properties:\\
(i) $A$ is mapped conformally onto $H$ by $F$, and 
$A$ contains a component $C = C(A)$ of $Y$ with no poles of 
$L$ on $\partial C$;\\
(ii) $C = C(A)$ contains a path $\gamma_C$ tending to infinity such that
\begin{equation}
{\rm Im} \, L(z) \geq |z|^{1/4} \quad
\hbox{as $z \to \infty$ on $\gamma_C$}
\label{H1}
\end{equation}
and
\begin{equation}
{\rm Im} \, z \to 0 \quad
\hbox{as $z \to \infty$ on $\gamma_C$.}
\label{H2}
\end{equation}
\end{lem}
\textit{Proof.}
As in \S\ref{Y} there exist $M$ distinct
components $K_j$, each contained in a component
$C_j$ of $Y$ which in turn lies in a component $A_j$ of
$W$. By Lemmas \ref{Wlem} and
\ref{Cjlem2} the number of $K_j$ contained in a given component $A$ of
$W$ is at most the number of $i$-points of $F$ in $A$, and 
this number is finite, while 
all but finitely many components $A$ of $W$ are
conformally equivalent to $H$ under $F$. Since
$M$ may be chosen arbitrarily large and $f$ has finitely many non-real
poles assertion (i) follows using Lemma \ref{lem1}. 

To prove assertion (ii) let $A$ and $C = C(A)$ be as in (i) and
observe that the function $u_C$ of 
(\ref{4}) is non-constant and subharmonic in the plane, but vanishes
on $\R$. Thus the existence of a path $\gamma_C$ satisfying (\ref{H1})
follows from a result of Barth, Brannan and Hayman \cite{BBH}.

It remains
to show that $\gamma_C$ also satisfies (\ref{H2}). To prove this
assume that the sequence $(w_\nu) \subseteq \gamma_C$ tends to infinity
with ${\rm Im} \, w_\nu \geq \varepsilon > 0$.
Take open discs $D_n$ of radius 
$\varepsilon$ about the poles $\zeta_n = (n + 1/2) \pi $ of
$\tan z $. Then $w_\nu \in \gamma_C \setminus \bigcup D_n$ and 
\begin{equation}
\tan w_\nu = O(1), \quad F(w_\nu) = \tan w_\nu + o(1) = O(1),
\label{H3}
\end{equation}
using (\ref{C3}) and (\ref{H1}). 
But $\gamma_C$ tends to infinity in $C \subseteq A$ and $F$ is
univalent on $A$, and so by passing to a subsequence
it may be assumed in view of (\ref{H3}) 
that $F(w_\nu) \to x \in \R$. Hence $\tan w_\nu \to x$,  
so that ${\rm Im} \, w_\nu \to 0$
using (\ref{hyp3}). This contradiction proves
(\ref{H2}).
\hfill$\Box$
\vspace{.1in}

Choose distinct components $A_1, \ldots , A_5$ as in Lemma \ref{lemG3}
and set $C_j = C(A_j)$, 
and take a large positive $R$ such that the circle $S(0, R)$
meets $\gamma_{C_j}$ for $j=1, \ldots , 5$. For each such $j$, choose
a subpath $\lambda_{C_j}$ of $\gamma_{C_j}$ lying in $|z| \geq R$
and joining $S(0, R)$ to infinity. It may
be assumed, after re-labelling if necessary,
that in $A^+(R, \infty)$ the path
$\lambda_{C_2}$ separates $\lambda_{C_1}$ from $\R$ and 
$\lambda_{C_3}$ separates $\lambda_{C_2}$ from $\R$ and, in view
of (\ref{H2}) and Lemma \ref{lemG3}(i), that
${\rm Im}\, z \to 0$ as $z \to \infty$ in $A_2 \cup A_3$. 
Denote positive constants by $c$, not necessarily the same at
each occurrence.

\begin{lem}\label{lemG4}
The function $F(z)$ satisfies the hyperbolic distance estimate
\begin{equation}
\frac{[ i, F(z)]_H }{|z|} \to \infty  
\quad
\hbox{as $z \to \infty$ in $A_2$.}
\label{H4}
\end{equation}
\end{lem}
\textit{Proof.} It follows from the fact that
${\rm Im}\, z \to 0$ as $z \to \infty$ in $A_2$
that for large $w$ in $A_2$
the largest disc of centre $w$ which lies in $A_2$ has radius $\tau (w) = o(1)$.
Choose $z^*$ in $A_2$ with $F(z^*) = i$ and
let $z \in A_2$ be large. Since the function
$u = G(w)$ in (\ref{uGdef})  
maps $A_2$ conformally onto $\Delta = D(0, 1)$, with inverse function
$w = h_1(u)$, the hyperbolic distance
$[0, G(z)]_\Delta$ is the infimum over all curves $\Gamma$ in $A_2$
joining $z^*$ to $z$ of
$$
\int_{G(\Gamma)} \frac2{1-|u|^2} \, |du| =
\int_\Gamma \frac{2 |dw| }{|h_1'(u)| (1 - |u|^2)} \geq
c \int_\Gamma \frac{|dw|}{\tau(w)} \geq \frac{|z|}{o(1)} ,
$$
using Koebe's theorem, and (\ref{H4}) follows since 
$[i, F(z)]_H = [0, G(z)]_\Delta$.
\hfill$\Box$
\vspace{.1in}

\begin{lem}\label{lemG5}
The component $C_2$ and its associated path $\lambda_{C_2}$
satisfy
\begin{equation}
{\rm Im}\, z \leq |z|^{-1/16} 
\quad
\hbox{as $z \to \infty$ on $\lambda_{C_2}$.}
\label{H5}
\end{equation}
\end{lem}
\textit{Proof.} Suppose that $z \in \lambda_{C_2}$ is large but
${\rm Im}\, z > |z|^{-1/16}$.
Then (\ref{C3}) and (\ref{H1}) give
\begin{equation}
| \tan z | \leq c |z|^{1/16}, 
\quad |F(z) - \tan z | \leq c |z|^{-1/8}, \quad
|F(z)| 
\leq c |z|^{1/16} .
\label{H6}
\end{equation}
Combining the last of these estimates with (\ref{hyp2}) and (\ref{H4}) forces
$$
\log \left( \frac1{{\rm Im}\, F(z) }  \right) \geq \frac{|z|}{o(1)} ,
\quad 
{\rm Im}\, F(z)  \leq \exp( - |z| ).
$$
Applying the second estimate of (\ref{H6}) again gives
$$
{\rm Im}\, ( \tan z ) \leq {\rm Im}\, F(z)  + c |z|^{-1/8}
\leq \exp( - |z| ) + c |z|^{-1/8} \leq c |z|^{-1/8} ,
$$
and using (\ref{hyp3}) this implies that
${\rm Im }\, z \leq c |z|^{-1/8}$,
contrary to assumption. This contradiction
proves the lemma.
\hfill$\Box$
\vspace{.1in}

Now let $u = u_{C_3}$ be the subharmonic function defined by (\ref{4}) and for
large $t > 0$ let $\theta (t)$ be the angular measure of $S(0, t) \cap
C_3$.
Then Lemma \ref{lemG5} and the choice of $C_3$ give 
$\theta (t) \leq c t^{-17/16} $
for large $t$, and 
$$
B(2r, u) \geq c \exp \left( \pi \int_{c}^{r} 
\frac{dt}{t \theta (t) } \right) \geq 
\exp \left( c r^{17/16} \right) 
$$
for large $r$. By (\ref{1}) and
(\ref{4}) this contradicts Lemma \ref{lemA2}, and 
the proof of Theorem \ref{thm1} is complete. 

\section{Concluding remarks}\label{contrast}

This section will outline some analogies and contrasts between the 
proofs of Theorems \ref{thmC} and \ref{thmA} (for $k=2$) on the one 
hand, and of Theorem \ref{thm1} 
on the other. In the proof of Theorem \ref{thm1} the function $F$ 
defined in (\ref{1}) plays a role comparable to that of the Newton
function $z - f(z)/f'(z)$ in
\cite{BEL,EdwH,SS}.
The connection between these apparently unrelated 
auxiliary functions may be seen as follows.
If $f$ is a real entire function such that $f$ and $f''+f$ 
have only real zeros then
following Frank's method \cite{Fra1,FH,FHP} write 
$$
f_1(z) = \cos z, \quad f_2 (z) = \sin z , \quad
W( f_1, f_2 ) = 1,
\quad
W(f_1, f_2, f) = f'' + f := \frac{f}{g^2} ,
$$
which gives, using standard properties of Wronskians
\cite[p.10]{Lai1},
$$
\frac1{(fg)^2} = W (f_1/f, f_2/f, 1) = W((f_1/f)',(f_2/f)' )
$$
and
$$
W(w_1, w_2) = 1,  \quad \hbox{where} \quad 
w_j = 
(fg) (f_j/f)' = f_j' g - L f_j g , \quad L = f'/f.
$$
It then follows that $w_1$ and $w_2$ are analytic in $H$ and that
the quotient $w_2/w_1$ has no critical points in $H$. But
$$
\frac{w_2}{w_1} = \frac{f_2' - Lf_2}{f_1'-Lf_1} =
\frac{f_1 - Lf_2}{-f_2-Lf_1} = \frac{TL-1}{L+T} = F ,
$$
where $T = \tan z $ and $F$ is as in (\ref{1}).
When $f$ and $f''$ have only real zeros
the same calculation with $f_1 (z) = 1, f_2 (z) = z$ leads to
$w_2(z)/w_1(z) = z - f(z)/f'(z)$, which is the Newton function.

The method in \S\ref{finiteorder} for the case of finite order
is closely related to the 
proof of Theorem \ref{thmC} \cite{EdwH,SS}, which as presented
for $f \in U_{2p}^*$ in \cite{EdwH} is lengthy.
It seems appropriate therefore to summarise the main steps for 
Theorem \ref{thmC} in the context of
the present method and to highlight the contrasts with Theorem \ref{thm1}. 
Indeed, suppose that $f \in U_{2p}^*$ with $2q$ distinct non-real zeros, and 
define $\phi$ and $\psi$ as in Lemma \ref{factorlem}. Then $\phi$
is a rational function \cite{EdwH}, 
and so $f' $ has finitely many non-real zeros, as has
$f''$ by repetition of the
same argument.
Let $L = f'/f$ as before and let $F$ be the Newton function of $f$.
Then $F$ has finitely many multiple points in $\C \setminus \R$
and finitely many non-real critical values, but in contrast to
the situation of Theorem \ref{thm1} the function $F$ has no 
asymptotic values in $\C \setminus \R$ (see e.g. \cite[Lemma 4]{Laedwhell}).
A standard argument
\cite[p.987]{BEL} then shows that $F$ is finite-valent on each component
$A$ of $W= \{ z \in H : F(z) \in H \}$. 
As in \S\ref{Y} there are at least $p+q$ components $K_j $ defined
as in Lemma \ref{componentslem}. 
Each $K_j$ satisfies $K_j \subseteq C_j \subseteq A_j$ 
for components $C_j$ of $Y = \{ z \in H : L(z) \in H\}$ and 
$A_j$ of $W $, not necessarily 
distinct. Moreover
each $K_j$ gives rise to a path $\Gamma_j \subseteq C_j$ as in \S\ref{Y}
on which $L(z) \to \infty$, and $\Gamma_j$ tends either to infinity
or to a non-real zero of $f$. Components $A_\nu$ of $W$ may then
be classified as type $(\alpha)$ or $(\beta)$ as in \S\ref{Y}, and
Lemma \ref{Cjlem3} applies to the type $(\beta)$ components.
Moreover if $\Gamma_j \to \infty$ then $F(z) \to \infty $ on
$\Gamma_j$ and so, since each $K_j$ has a zero of $L$ and so a
pole of $F$ on its boundary, the valency of $F$ on a type $(\alpha)$
component $A$ of $W$ exceeds the number $\mu_A$ of $K_j$ contained
in $A$ by at least $1$, so that by the Riemann-Hurwitz formula
the number of critical points of
$F$ in $A$ is at least $\mu_A$. Since these critical points are either
non-real zeros of $f$ or of $f''$,
it follows that there are at least $p$ zeros of $f''$ in $H$.

\vspace{.1in}
\noindent
School of Mathematical Sciences,
University of Nottingham,
NG7 2RD.\\
jkl@maths.nottingham.ac.uk

\end{document}